\documentclass[journal]{IEEEtran}

\usepackage{pgfplots,pgfplotstable,booktabs}

\pgfplotsset{compat=1.5}

\usepackage{amsmath, amssymb}

\usepackage{caption}

\usepackage{graphicx}

\usepackage{subfigure}

\usepackage{graphicx,epstopdf}
\DeclareGraphicsExtensions{.eps}

\usepackage{etoolbox}

\apptocmd{\sloppy}{\hbadness 10000\relax}{}{}

\apptocmd{\sloppy}{\hbadness 6236\relax}{}{}

\usepackage{cite}

\usepackage{epsf} 

\usepackage{tikz}
\newcommand*\circled[1]{\tikz[baseline=(char.base)]{
            \node[shape=circle,draw,inner sep=1pt] (char) {#1};}}

\begin{document}

\title{Voltage Multistability and Pulse Emergency Control for Distribution System with Power Flow Reversal}

\author{Hung~D.~Nguyen,~\IEEEmembership{Student Member,~IEEE,} and~Konstantin~Turitsyn,~\IEEEmembership{Member,~IEEE}% <-this % stops a space
\thanks{Hung D. Nguyen and Konstantin Turitsyn are with the Department of Mechanical Engineering, Massachusetts Institute of Technology, Cambridge, MA, 02139 USA e-mail: hunghtd@mit.edu and turitsyn@mit.edu}}% <-this % stops a space

 \markboth{IEEE Transactions on Smart Grid 2015, in press}%
 {Shell \MakeLowercase{\textit{et al.}}: Voltage Multistability and Pulse Emergency Control for Distribution System with Power Flow Reversal}%

\maketitle

\begin{abstract}
%\boldmath
High levels of penetration of distributed generation and aggressive reactive power compensation may result in the reversal of power flows in future distribution grids. The voltage stability of these operating conditions may be very different from the more traditional power consumption regime. This paper focused on demonstration of multistability phenomenon in radial distribution systems with reversed power flow, where multiple stable equilibria co-exist for the given set of parameters. The system may experience transitions between different equilibria after being subjected to disturbances such as short-term losses of distributed generation or transient faults. Convergence to an undesirable equilibrium places the system in an emergency or \textit{in extremis} state. Traditional emergency control schemes are not capable of restoring the system if it gets entrapped in one of the low voltage equilibria. Moreover, undervoltage load shedding may have a reverse action on the system and can induce voltage collapse. We propose a novel pulse emergency control strategy that restores the system to the normal state without any interruption of power delivery. The results are validated with dynamic simulations of IEEE $13$-bus feeder performed with SystemModeler software. The dynamic models can be also used for characterization of the solution branches via a novel approach so-called the admittance homotopy power flow method.

\end{abstract}


 \maketitle

%\graphicspath{{}}

\section{Introduction}

The increasing levels of penetration of distributed generators (DGs), either renewable or gas-fired will cause the distribution grids to operate in unconventional conditions. The flow of active or reactive power may become reversed in certain realistic situations such as sunny weekday time in residential areas with high penetration of photovoltaic panels. Active participation of future distribution level power electronics in reactive power compensation may also lead to the local reversal of reactive power flows. These kind of operating conditions are not common to existing power grids, but may become more common in the future and may also have a serious effect on the overall voltage stability of the system.

The strong nonlinearities present in the power system determine the existence, multiplicity, and stability of the viable operating points \cite{Hiskens95,Hiskens89}. The nonlinear control loops inside individual system components are responsible for the voltage collapse \cite{dobson1992voltage,dobson92collapse,Cutsem} and loss of synchrony phenomena \cite{BialekBook} that have caused some of the most severe blackouts in the recent history. %\cite{BlackoutPapersbyCanizaresDobson.etal}. 
Generally, the power flow equations that are commonly used for the description of steady states of the power system \cite{Kundur} may have multiple solutions \cite{Thorp}, but in typical operating conditions, there always exists a high voltage solution that is considered a normal operating point \cite{Chiang90}. 

The power flow equations solutions manifold has been studied rather extensively in the context of transmission grids; nevertheless, the structure of the solution manifold in distribution grids in reversed power flow regime is however poorly understood, although there are reasons to believe that it will be very different from the classical nose-curve type manifold. Even though the direction of the power flow does not affect the qualitative properties of the solutions in linear (DC power flow) approximation, it becomes important when the nonlinearity is strong. The symmetry between the normal and reversed power flow solutions is broken because the losses that are the major cause of nonlinearity in the power flow equations are always positive. In traditional distribution grids the consumption of power and the losses have the same sign, while in the situation with reversed flows the processes of power injection and thermal losses are competing with each other. This competition may manifest itself in the appearance of new solutions of power flow equations that do not exist in the non-reversed power flow regime. From power engineering perspective, this phenomena can be understood with the following argument. In the presence of power flow reversal, the power injections raise the voltage to high enough levels for low voltage equilibria to appear. The existence of low voltage equilibria may be demonstrated by continuation type rigorous mathematical arguments, that we have presented in Appendix \ref{app:newsol}. It is also based on our observations from numerical simulations/experiments and the discussions of the existence of power flow problem presented in \cite{thorp1986reactive,ilic1992network}. This phenomena was observed by one of the authors in a recent work \cite{wang2012distflow} but has not been explored in greater details since then.

Even for the traditional nose-curve scenario, the second low voltage solution may be stable under some conditions. This has been recognized for a long time \cite{Overbye1994, Nose_Taranto, Venkatasubramanian1992, hill1994stability}. Moreover, Venkatasubramanian et al. in \cite{Venkatasubramanian1992} noted that the situations in which the systems gets trapped at the second stable equilibrium have been observed. However, the relevance of the low voltage stable equilibrium did not draw much attention and/or has not been studied extensively because this stable equilibrium is neither viable nor convincingly verified numerically due to modeling difficulties. The main problem in the assessment of the stability is the highly complex nature of the load dynamics. The dynamic behavior of the loads are a result nonlinear interactions of millions of heterogeneous components that are poorly understood and not fully known to the operator of any given grid. At the same time, the dynamic behavior has a direct effect on the stability properties that cannot be directly assessed via static power flow analysis \cite{Overbye1994}. In this work, we address this problem by studying the stability with dynamic load models that are consistent with existing models in normal conditions but does not suffer from the convergence problems in abnormal situations.

Although, new equilibria are not suited to the normal operation of the power system, they may cause serious effects on the transient stability and post-fault recovery of the system. The effect is similar although more serious than the Fault Induced Delayed Voltage Recovery (FIDVR) observed in power grids with high share of induction motors in the load composition \cite{Hiskens2010,Hiskens2005,Dobson2012}. Stalling of induction motors may cause the delay in the restoration process. At the same time, when the grid has stable low voltage equilibria, the system may get entrapped and fail to escape from the equilibrium at all. In this case, the likely outcome of the dynamics will be the tripping action of the undervoltage protective relays and consequent partial outage of the power grid. Hence, it is important to revisit the voltage protection controls for the future power grids with high penetration of distributed generation. In this work we proposed specific pulse emergency control strategies (PECS) that are designed to restore the system to the normal operating condition.

The key contributions of this work are summarized below. In section \ref{sec: Dynamic Simulations 3 bus}, we perform dynamic simulations on a three-bus model that illustrate the effect of multistability and possibility of system entrapment in undesirable low voltage states. The simulations illustrate both the stability of some regions of low voltage part of the nose curve as well as the effect of multistability at high load levels. Next, in section \ref{sec: Dynamic Simulations 13 bus}, we introduce a novel admittance homotopy power flow technique to find multiple solutions to the power flow equations, and perform dynamic simulations of a larger system, the IEEE thirteen-bus test feeder. We dedicate section \ref{sec:emergencycontrol} to revisit the current emergency control actions and design the pulse emergency control for multistability. The importance of proper power reversal regulations as well as the proper assessment of DG penetration level on planning stages is also discussed in the end of section \ref{sec:emergencycontrol}.

\section{Dynamic simulations of a three-bus network} 
\label{sec: Dynamic Simulations 3 bus}
In this and the next section, we introduce the model and perform dynamic simulations of two radial networks to show that the two stable equilibria of load dynamics equations may coexist at the same time, and that the distribution system may become entrapped at the lower voltage equilibrium. 

\begin{figure}[ht]
    \centering
    \includegraphics[width=0.8 \columnwidth]{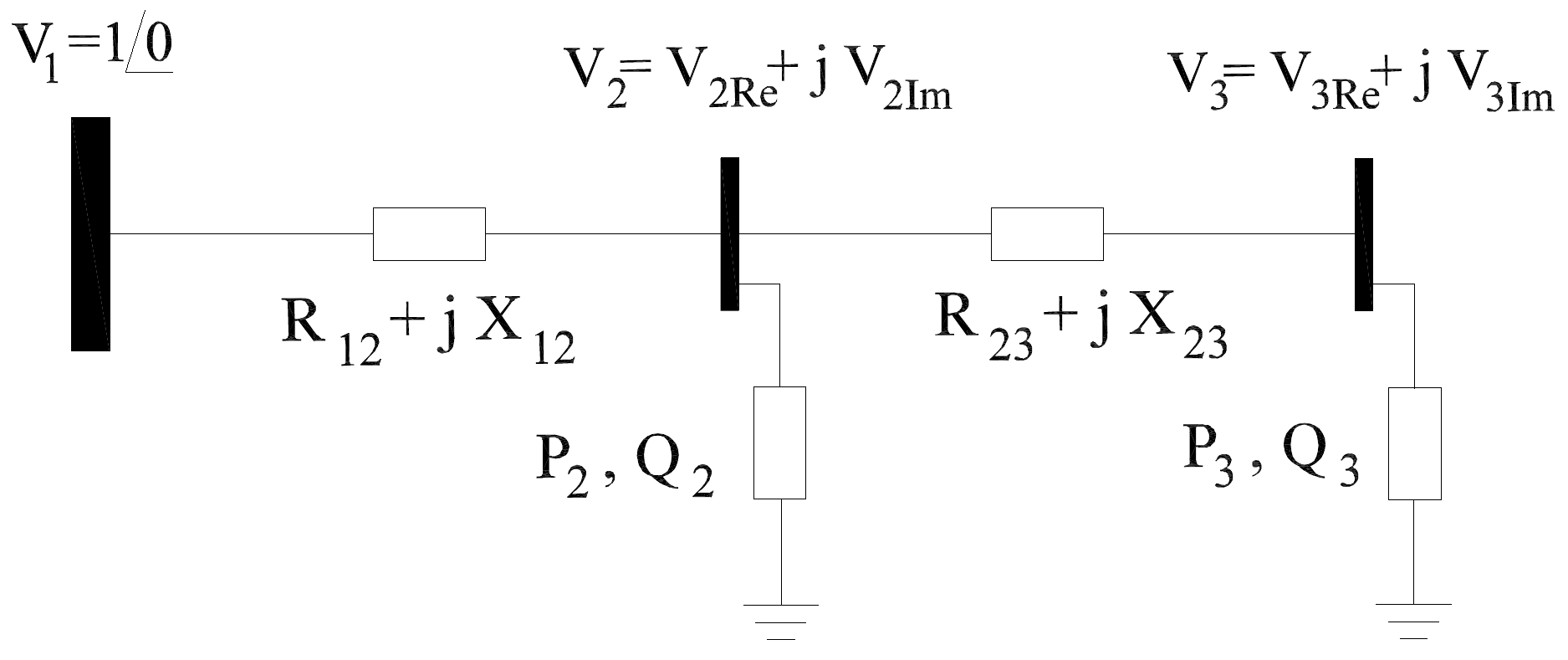}
	\caption{A three-bus network}
    \label{3bus}
\end{figure}

First, we consider a three-bus network as shown in Figure \ref{3bus} with bus $1$ being the slack bus and buses $2$ and $3$ representing the dynamic loads with distributed generation exporting reactive power. This system could represent the future distribution grids with the inverters of PV panels participating in voltage regulation (see \cite{turitsyn2011options} for further discussions of this proposal). Alternatively, it could represent a highly capacitive grid, for example involving long underground cables. 

The loads are modeled as dynamic admittances with conductance $g$ and susceptance $b$ that evolve according to the following dynamic equations:
\begin{equation} \label{eq:gg}
\centering
\tau_1\dot{g}=-(p-P^s)
\end{equation}

\begin{equation} \label{eq:bb}
\centering
\tau_2\dot{b}=-(q-Q^s)
\end{equation}

The values $P^s$ and $Q^s$ describe the static (steady-state) power characteristics of the loads that are achieved in equilibrium. This form of dynamic load model is consistent with the most common ones in normal operating conditions, but is also applicable to highly nonlinear transients. The traditional load models \cite{Karrison94, Hill93} with time-varying power consumption levels may not have solutions in all the transient states, whereas the model based on local admittance levels always has solutions.

We assume that the steady-state active and reactive power levels do not depend on the bus voltage, so that $P^s$ and $Q^s$ depend only on time. In other words, the loads can be classified as constant power loads, that attempt to achieve the given power demand levels $P^s(t), Q^s(t)$. This assumption is clearly a simplification of the real loads. However, it is a common assumption in most of the classical voltage stability studies, and also may be a good approximation of power systems with aggressive VAR compensation or power grids interconnected through fast voltage source converters (VSC). Moreover, the qualitative nature of the results does not depend on the specific static characteristic and more realistic models are used in next section.

In our simulations, we use $\tau_1 = 3 \,s$, $\tau_2 = 0.001\, s$ for the load at bus $ 2$ and $\tau_1 = \tau_2 = 0.01\, s$ for the load at bus $3$. The actual values of time constants do not affect any stability properties, we chose them for convenience of presentation, but physically the fast dynamic loads could correspond to power electronics regulating the voltage levels on the consumer side. The network parameters are given as the following. The impedance of the line $1-2$, $z_{12} = 0.03 + j 0.15 \,p.u.$; the impedance of the line $2-3$, $z_{23} = 0.33 + j 1.65\,p.u.$. The load power consumption levels are defined as: $P_i=P_{G\,i}+P_{L\,i}$ and $Q_i=Q_{G\,i}+Q_{L\,i}$; where $i$ is the load number, $i=2,3$. $P_{G\,i}$ and $Q_{G\,i}$ are the active and reactive powers produced from distributed generators at bus $i$, whereas, $P_{L\,i}$ and $Q_{L\,i}$ are the active and reactive powers consumed at bus $i$. The power levels have negative values if the bus is generating power, whereas, positive values indicate that the bus is consuming power. Hereafter, all simulations rely on the assumption that all components/devices are able to operate through short periods under low voltage conditions, and without disconnection to the grid. In other words, all components/devices have enough Low Voltage Ride Through (LVRT) capabilities. This assumption holds when high levels of system reliability are enforced. In this case, we assume that DGs continue to export power in the post-fault recovery.

In this paper, we mainly focus on non-synchronous DGs, i.e. wind turbines and solar panels. These DGs can be modeled as either $PQ$ or $PV$ bus by selecting corresponding control mode \cite{Erlich2006, Tamimi2013}. As a common practice, DGs can be switched from constant voltage to constant reactive power output when reactive power limits are reached \cite{Tamimi2013}. We use the $PQ$ model or in equivalent, a negative load equipped with integral controller ($I$), $K_I = \tau$. The aggregated loads represent the combination of both the traditional consumption loads and DGs represented as negative loads. For the sake of simplicity, we assign a single time constant for the aggregated one. More sophisticated voltage control strategies may be explored in future more detailed studies. Even though the existing standards for DGs integration do not allow them to control voltage, i.e. require them to maintain a unity power factor; in the future, it is expected to change.

\subsection{Transient from high voltage equilibrium to the low voltage one} \label{H2L}
The objective of our first scenario simulations is to demonstrate that the system may get entrapped in the stable low voltage equilibrium.

\begin{figure}[ht]
    \centering
    \includegraphics[width=1 \columnwidth]{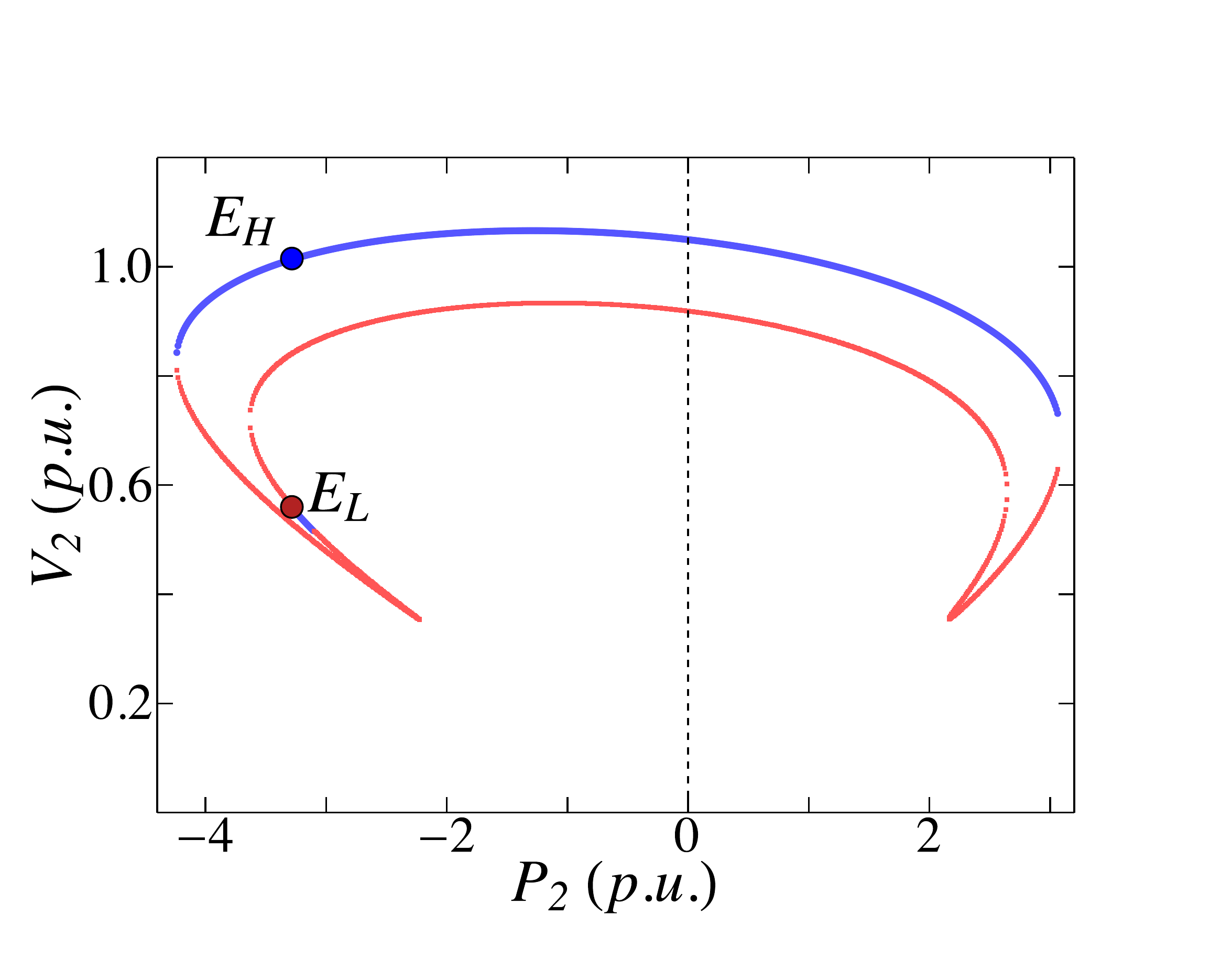}
	\caption{Voltage multi-stability in a three-bus network}
    \label{fig:voltagemultistability}
\end{figure}

The bus $3$ has base demand level $P_3^s = -0.189 \,p.u.$ and $Q_3^s = -0.222\, p.u.$ which corresponds to a capacitive load producing active power. The base active and reactive power demands of load $2$ are also given as $P^s_2= -3.284\, p.u.$, $Q_2^s = -0.167\, p.u.$. While keeping $Q_2^s \,,P_3^s\,,Q_3^s$ fixed to equal to the base level and changing the active demand level at bus $ 2$, $P_2^s$, we can plot different $PV$ curves as in Figure \ref{fig:voltagemultistability} with two stable equilibria. In these plots, the blue dot segments represent stable equilibria and the red dot segments - the unstable ones, all observed for different values of $P_2^s$. The large blue dot represents the high voltage stable equilibrium, $E_H$, and the large red one marks the low voltage one, $E_L$. For the given parameters, the two stable equilibria, $E_H$ and $E_L$, correspond to two levels of voltage at bus $2$, i.e. $V_2 = 1.012\, p.u.$ and $V_2 = 0.560 \,p.u.$, respectively. The following scenario initiates the transition of the system from the high voltage equilibrium to the low voltage one.  

\begin{figure}[!ht]
    \centering
    \includegraphics[width=.9 \columnwidth]{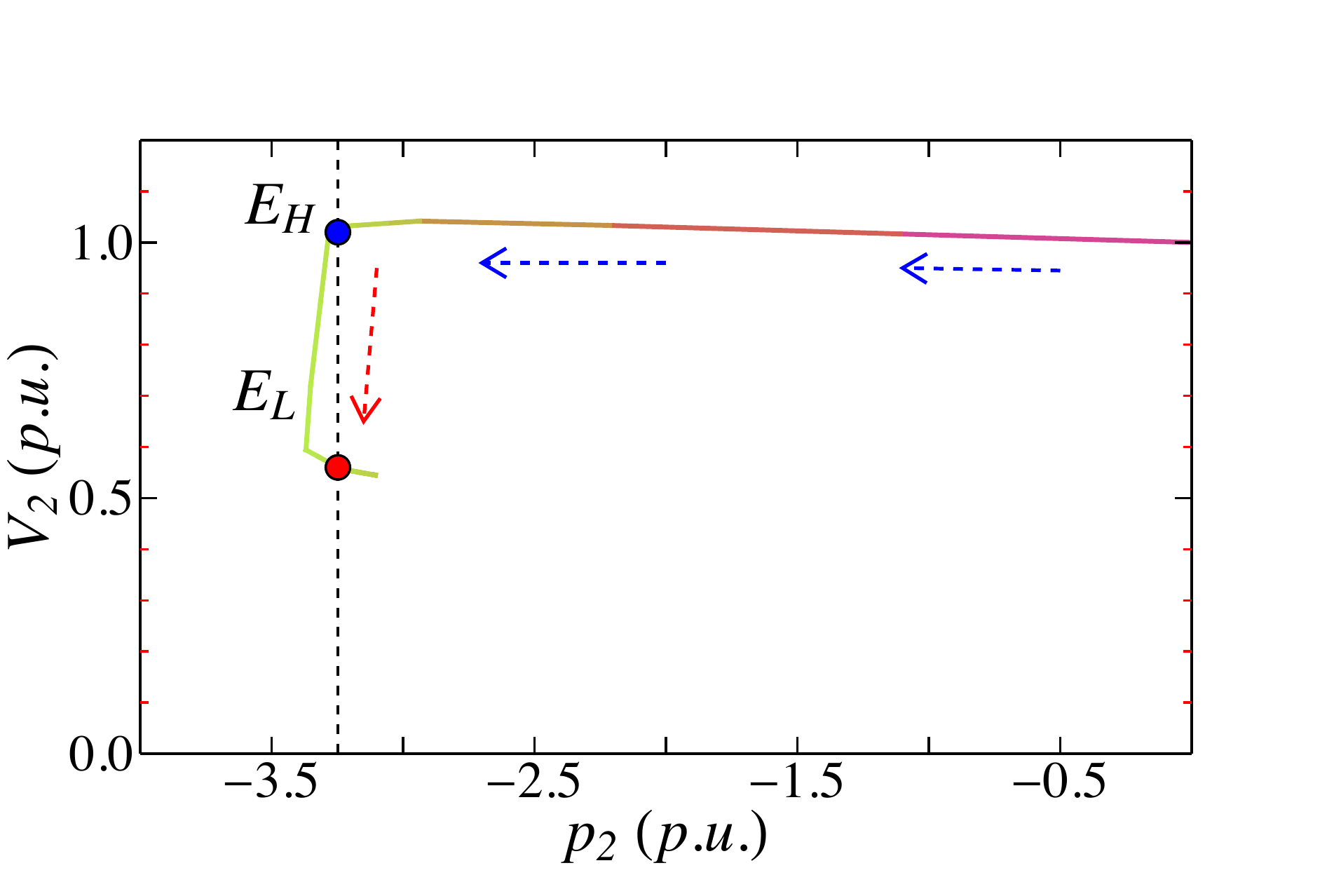}
	\caption{The $PV$ curve at bus $2$, $t<25\, s$}
    \label{fig:p2v24000s_H2L}
\end{figure}

Initially the system is operating at the high voltage equilibrium, $E_H$. In the dynamic simulation, the preferred operating condition can be reached by choosing appropriate initial conditions in the neighborhood of the steady state. In this case, zero initial condition is suggested. The initial transient to the equilibrium point is shown in Figure \ref{fig:p2v24000s_H2L} following the blue arrows. 

After the system reaches the high voltage equilibrium, a large disturbance occurs, i.e. distributed generation is partially lost, at $t_{d1}=15\,s$. For example it could represent the cloud covering the PV panels with a shadow or action of protective relays after a short-circuit. The aggregate load on bus $2$ changes its behavior from generation to consumption of active and reactive power. As a result, the system starts to move away from the high voltage stable equilibrium, $E_H$, and approach the low voltage one, $E_L$. This process is presented in Figure \ref{fig:p2v24000s_H2L} with red arrows. The transient dynamics of the system dies out around $t\approx{15.3}\, s$. 

\begin{figure}[ht]
    \centering
    \includegraphics[width=1 \columnwidth]{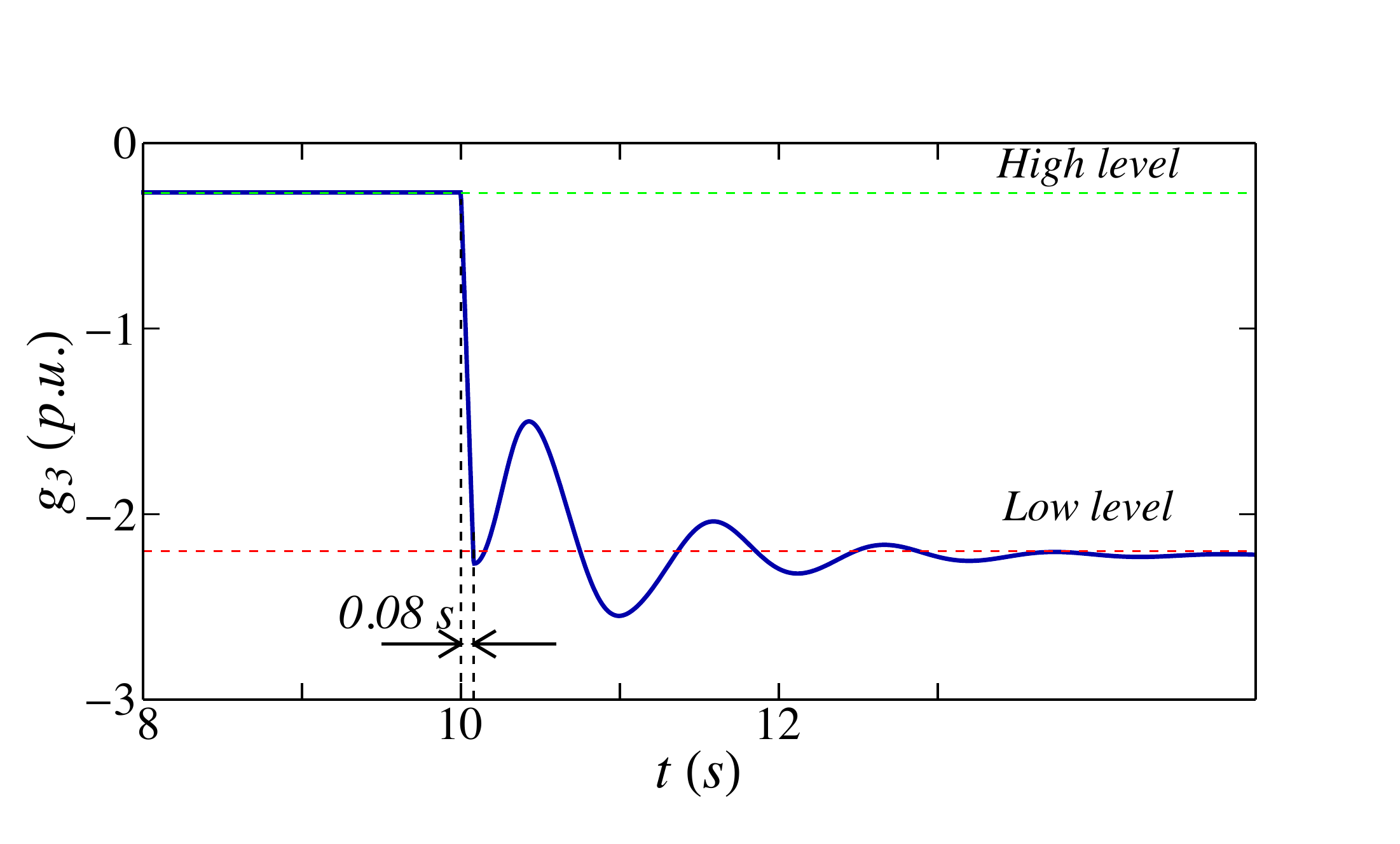}
	\caption{Conductance dynamics at bus $3$ due to a transient short circuit}
    \label{fig:shortcircuit}
\end{figure}

The same transient from the high voltage equilibrium to the low voltage equilibrium can be also observed in another system where $z_{12} = z_{23} = 0.1464 + j0.5160\, p.u.$, $P^s_2=-0.7\, p.u.$, $Q_2^s = -0.9\, p.u.$, $P_3^s =- 0.75 \,p.u.$, $Q_3^s = -0.45\, p.u.$; $\tau_1 = \tau_2 = 0.07\, s$ for the load at bus $ 2$ and $\tau_1 = \tau_2 = 0.03\, s$ for the load at bus $3$; if a transient fault located at either bus $2$ or bus $3$ occurs at $t_{d1}=10\,s$ for a short time, e.g. $0.08\,s$. The conductance of the load at bus $3$ during the transient is shown in Figure \ref{fig:shortcircuit}. The high level and low level refer to the voltage level of the corresponding equilibrium.

\begin{figure}[!ht]
    \centering
    \includegraphics[width=0.8 \columnwidth]{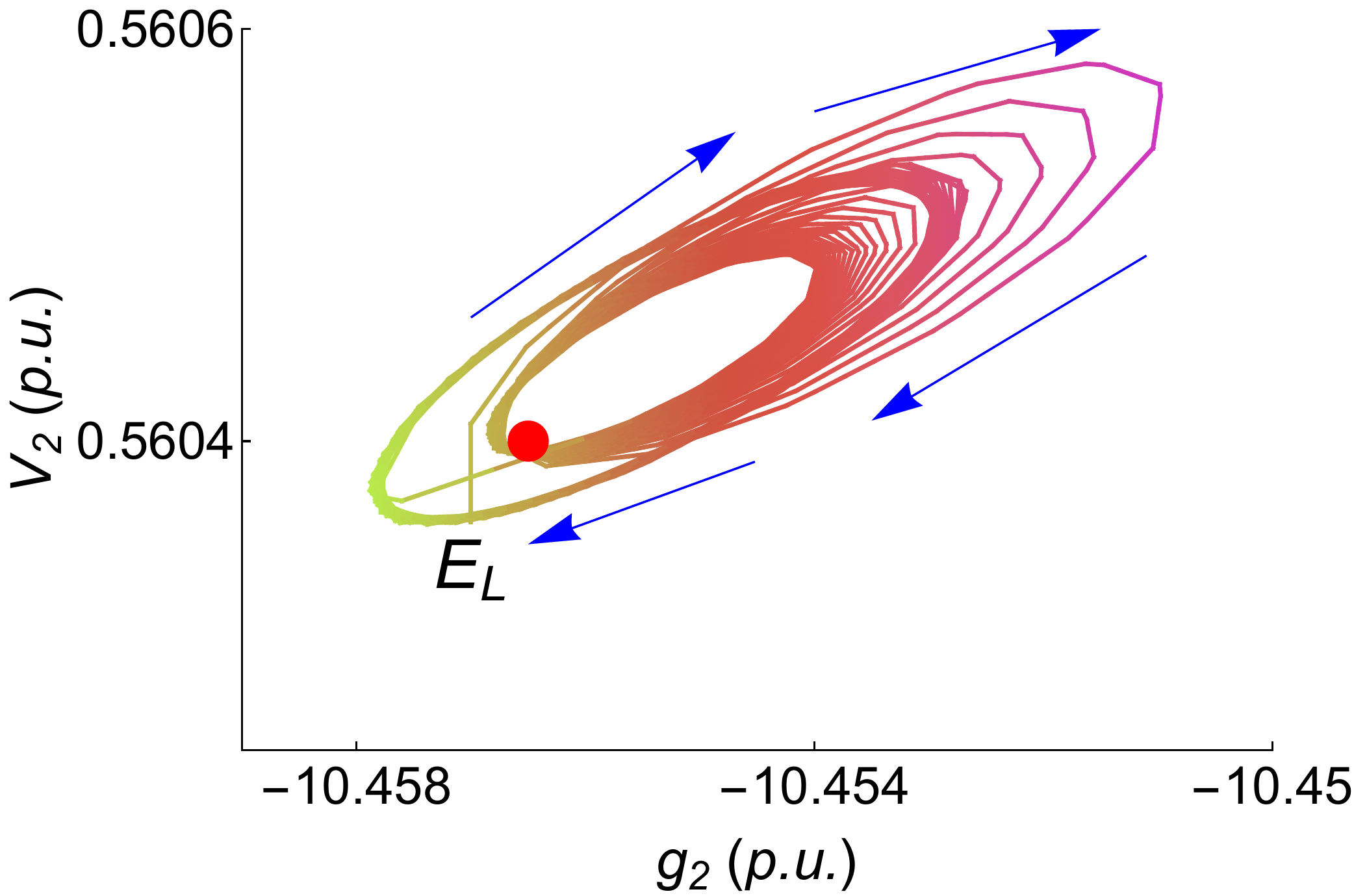}
	\caption{The $gV$ curve at bus $2$, $25\,s\leq t\leq{30\, s}$}
    \label{fig:p2v2ts_H2L}
\end{figure}

Later, at $t_{d2} = 25\, s$ the second disturbance in $P_2^s$ occurs that changes the demand $P_2^s$ to some higher value $P_2^s = -3.286\, p.u.$ for $0.01\,s$. As shown in Figure \ref{fig:p2v2ts_H2L}, the system first moves away from the low voltage equilibrium following the blue arrows then returns back to the same equilibrium following the blue dashed arrows. Our numerical experiments with several disturbances with different amplitudes and durations indicate that the low voltage equilibrium is indeed non-linearly stable and is characterized by the region of attraction with finite size.

\section{Dynamic simulations of the IEEE thirteen-bus distribution feeder test case} \label{sec: Dynamic Simulations 13 bus}
\begin{figure}[!ht]
    \centering
    \includegraphics[width=0.8 \columnwidth]{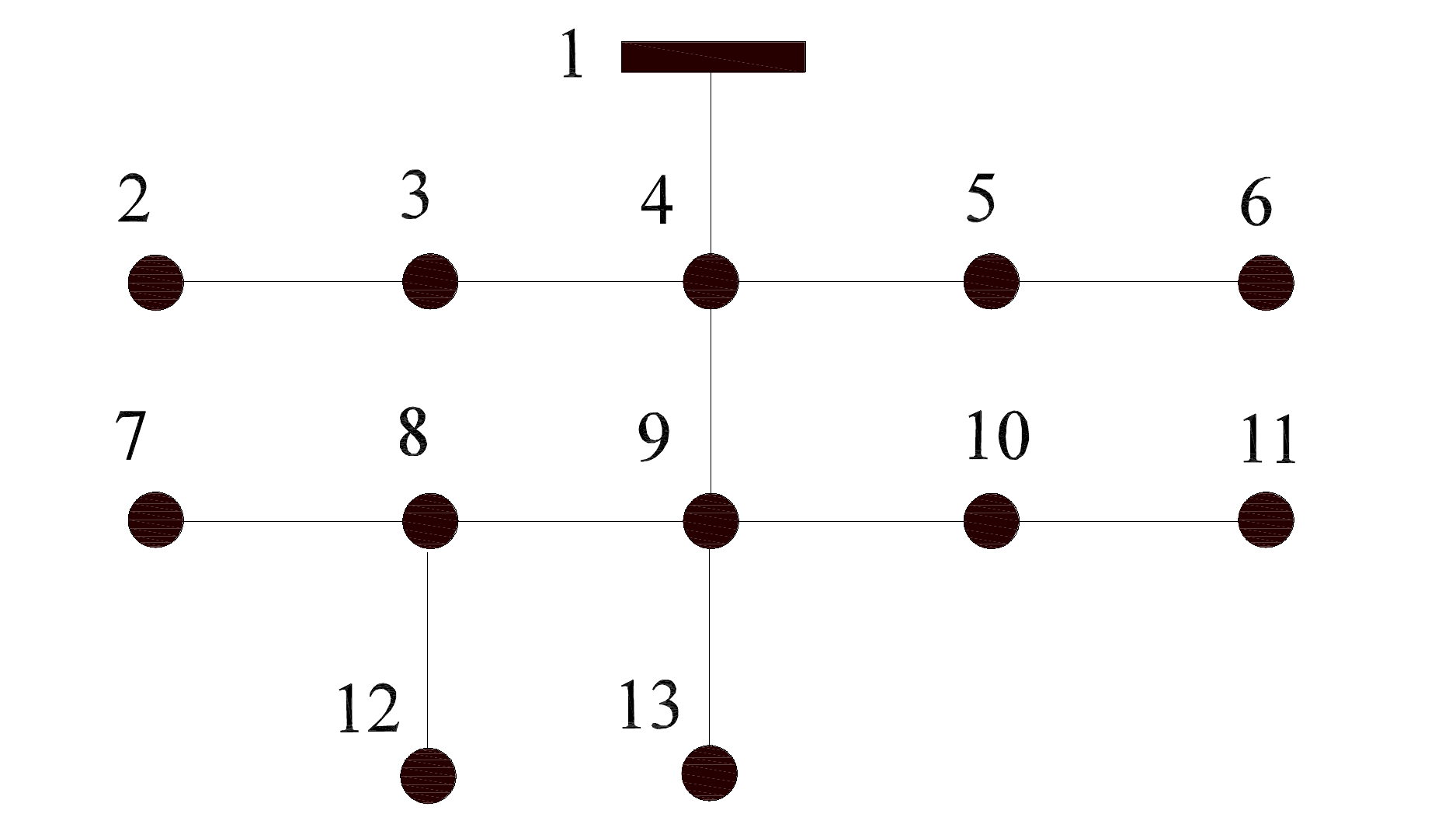}
	\caption{The thirteen-bus feeder test case}
    \label{fig:13busfeeder}
\end{figure}

In this section, we study the multistability phenomenon and the transition among equilibria in a more realistic thirteen-bus radial system based on the IEEE thirteen-bus distribution feeder test case described in \cite{IEEEtestcase}. The network configuration is shown in Figure \ref{fig:13busfeeder} and the branch data can be found in Table \ref{branchdata}.

\begin{table}[ht]
    \caption{The thirteen-bus test feeder branch data}
    \centering
    \label{branchdata}
    \begin{tabular}{|c |c |c||c |c |c|}
    \hline
    \textbf{From} & \textbf{To} & \textbf{Line impedance}& \textbf{From} & \textbf{To} & \textbf{Line impedance} \\
    \textbf{bus} & \textbf{bus } & \textbf{$(p.u.)$}&  \textbf{bus } & \textbf{bus} & \textbf{$(p.u.)$}\\
               \hline
    $1$ & $4$ & $0.041+j0.131$&$8$ & $9$ & $0.027+j0.070$ \\
    $2$ & $3$ & $0.027+j0.070$ &$4$ & $9$ & $0.041+j0.131$  \\
    $3$ & $4$ & $0.045+j0.100$ &$9$ & $10$ & $0.006+j0.020$\\
    $4$ & $5$ & $0.032+j0.070$ &$10$ & $11$ & $0.027+j0.070$\\
    $5$ & $6$ & $0.100+j0.240$ &$8$ & $12$ & $0.039+j0.090$\\
    $7$ & $8$ & $0.015+j0.050$ &$9$ & $13$ & $0.020+j0.065$\\
    
    \hline
    \end{tabular}
\end{table}

In the test case, bus $1$ is the slack bus and bus $2$ is a dynamic load with $\tau_1=\tau_2=0.01 \,s$. The other buses are modeled as polynomial loads with load bus $k$ is described as:
\begin{eqnarray} \label{zip1}
%\begin{split} 
P_k=P_{k\,0}(a_{k\,P}|V_k|^2+b_{k\,P}|V_k|+c_{k\,P}) \\
Q_k=Q_{k\,0}(a_{k\,Q}|V_k|^2+b_{k\,Q}|V_k|+c_{k\,Q}) 
%\end{split}
\end{eqnarray}
where $a_{k\,P}+b_{k\,P}+c_{k\,P}=a_{k\,Q}+b_{k\,Q}+c_{k\,Q}=1$, and $P_{k\,0}$ and $Q_{k\,0}$ are the active power and reactive powers the load $k$ consumes at the reference voltage level, $V_0=1\,p.u.$ \cite{Cutsem}. The load consumption levels are listed as $P_{2\,0}=-0.85$, $Q_{2\,0}=0.1$, $P_{3\,0}=P_{4\,0}=P_{5\,0}=P_{6\,0}=0.1$, $Q_{3\,0}=Q_{4\,0}=Q_{5\,0}=Q_{6\,0}=0.1$, $P_{7\,0}=Q_{7\,0}=0.5$, $P_{8\,0}=0.5$, $Q_{8\,0}=0.3$, $P_{9\,0}=0$, $Q_{9\,0}=-1$, $P_{10\,0}=0$, $Q_{10\,0}=-1$, $P_{11\,0}=Q_{11\,0}=0$, $P_{12\,0}=0$, $Q_{12\,0}=-0.5$, $P_{13\,0}=0$, $Q_{13\,0}=-1$. Simply, let $a_{k\,P}=b_{k\,P}=0.01$, $a_{k\,Q}=b_{k\,Q}=0.005$, and $c_{k\,P}=c_{k\,Q}=0.985$, $k=\overline{3,13}$. Hence, bus $2$ exports active power and consumes reactive power. Bus $9$, $10$, $12$, and $13$ consume active power and export reactive power. The other buses consume both active and reactive power. To demonstrate the phenomenon of multistability, we first need to find multiple solutions of the network for a given set of parameters which are, in this case, the load consumption levels. 

\subsection{Admittance homotopy power flow} \label{sec:dynamicpowerflow}
To find multiple solutions of the power flow problem in large scale grids we introduce a novel homotopy type technique that is described in details in this section. 

Within the framework of our approach, we represent each load or generation with an unknown effective admittance $y_k = g_k + i b_k$ and introduce the standard power flow constraints
\begin{eqnarray} \label{pc}
& P_k=g_k|V_k|^2 \\ 
& Q_k=b_k|V_k|^2 \label{qc}
\end{eqnarray}
with $P_k, Q_k$ fixed at the desirable consumption levels and the voltage satisfying the standard Kirchhoff current law of the form:

\begin{equation} \label{eq:KCL}
\sum_{l=1}^ny_{kl}V_l=y_kV_k
\end{equation}

where $y_{kl}$ is the entry of the admittance matrix $Y_{bus}$ that corresponds to the link between bus $k$ and bus $l$; $V_k$ and $V_l$ are the voltage at bus $k$ and bus $l$, respectively. Bus $1$ is the slack bus of the considered $n$-bus system. \eqref{eq:KCL} can be rewritten as:
\begin{equation*}
\sum_{l=1}^n[y_{kl}V_l-y_k\delta_{kl}]V_l=0
\end{equation*}

Or:

\begin{equation} \label{eq:newPL}
Y_{kl}(\underline{y}) V_l = -Y_{k1} V_1
\end{equation}
where $\delta_{kl}=1$ if $k=l$ and $\delta_{kl}=0$ otherwise.

\eqref{eq:newPL} defines a system of algebraic equations for the vector of effective admittances $\underline{y}$ whose solution represents one of the branches of power flow equation solutions. Note, that this system equation is completely equivalent to the standard load flow equations, it is just written in terms of different variables, or in other words, the solution manifold is parameterized with different set of coordinates.

To characterize all the branches of the solution we fix all of the power flow constraints except for one, say $P_2 = g_2 |V_2|^2$ and begin to vary the value of $g_2$ thus exploring one of the cross-sections of the solution manifolds. For almost every value of $g_2$, we can find the solution of the nonlinear power flow constraints (\ref{pc},\ref{qc}), and calculate the corresponding value of $P_2$. Thus, this way we define the parametric representation of the solution manifold cross-section. Like in the traditional homotopy approaches, there will be some points at which the Jacobian of the (\ref{pc},\ref{qc}) with respect to the free parameter $g_2$ becomes singular where one needs some reconditioning or different initial conditions to jump on a different branch. However, the key advantage of the admittance based homotopy is that this singularity occurs at low voltage levels, so the method can characterize the standard nose parts of the curves without any divergence problems.

\begin{figure}[!ht]
    \centering
    \includegraphics[width=1 \columnwidth]{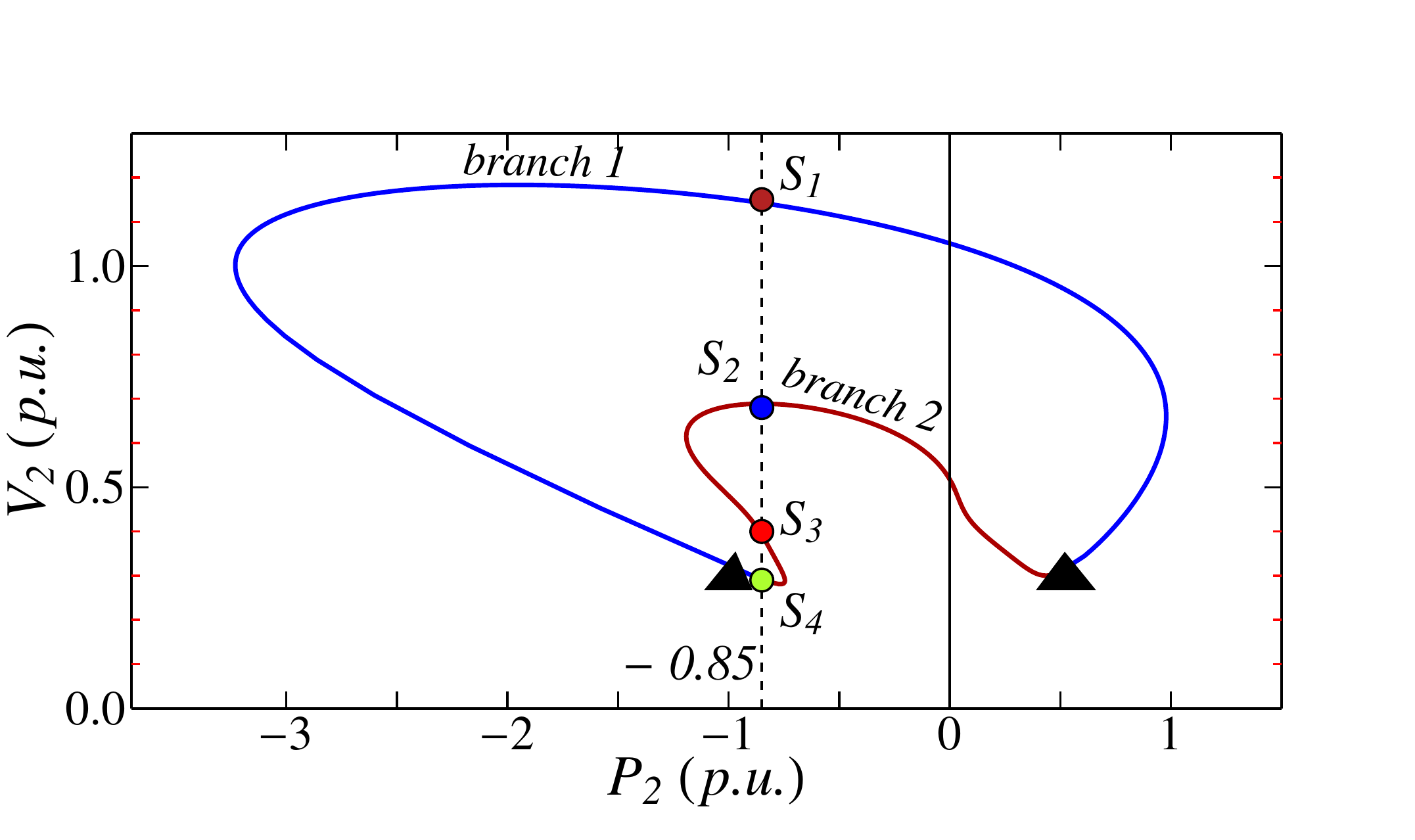}
	\caption{The thirteen-bus feeder test case nose curve}
    \label{fig:13busnose}
\end{figure}

In Figure \ref{fig:13busnose}, we show all the branches of the power flow solutions for the $13$-bus system obtained using the  admittance homotopy approach.

The two branches of the nose curve terminate at the singular points which are labeled with ``$\blacktriangle$" symbols. After the set of nose curves $P_2-V_k$, $k=\overline{2,13}$ are obtained, one can determine the solutions to the power flow problem. For example, from Figure \ref{fig:13busnose}, one can find that there are four solutions, $S_1, S_2, S_3, S_4$ corresponding to four voltage levels, $V_2=1.15,\,0.69,\,0.4,\,0.3\,(p.u.)$, respectively. The transition among these equilibria is simulated in the next subsection.

The proposed admittance homotopy power flow method does not suffer from divergence problems at the maximal loading points as traditional solution approaches based on iterative techniques such as Newton-Raphson method. Also, unlike algebraic geometry approaches, the dynamic homotopy power flow method can be scalable to large scale systems.

\subsection{Transition among equilibria}
\begin{figure}[!htb]
    \centering
    \subfigure[]{{\includegraphics[width=4.5cm]{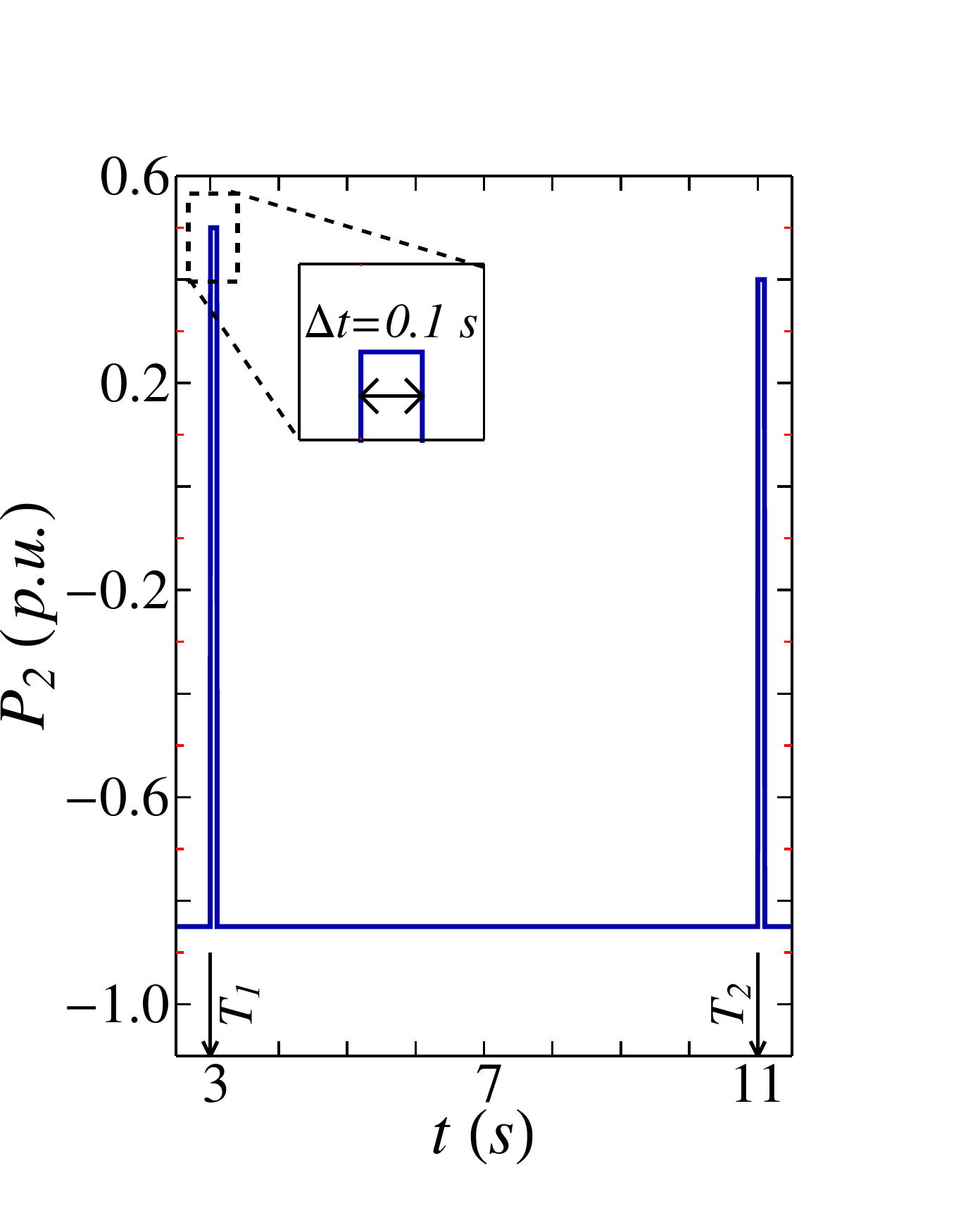}}}%
    \subfigure[]{{\includegraphics[width=4.3cm]{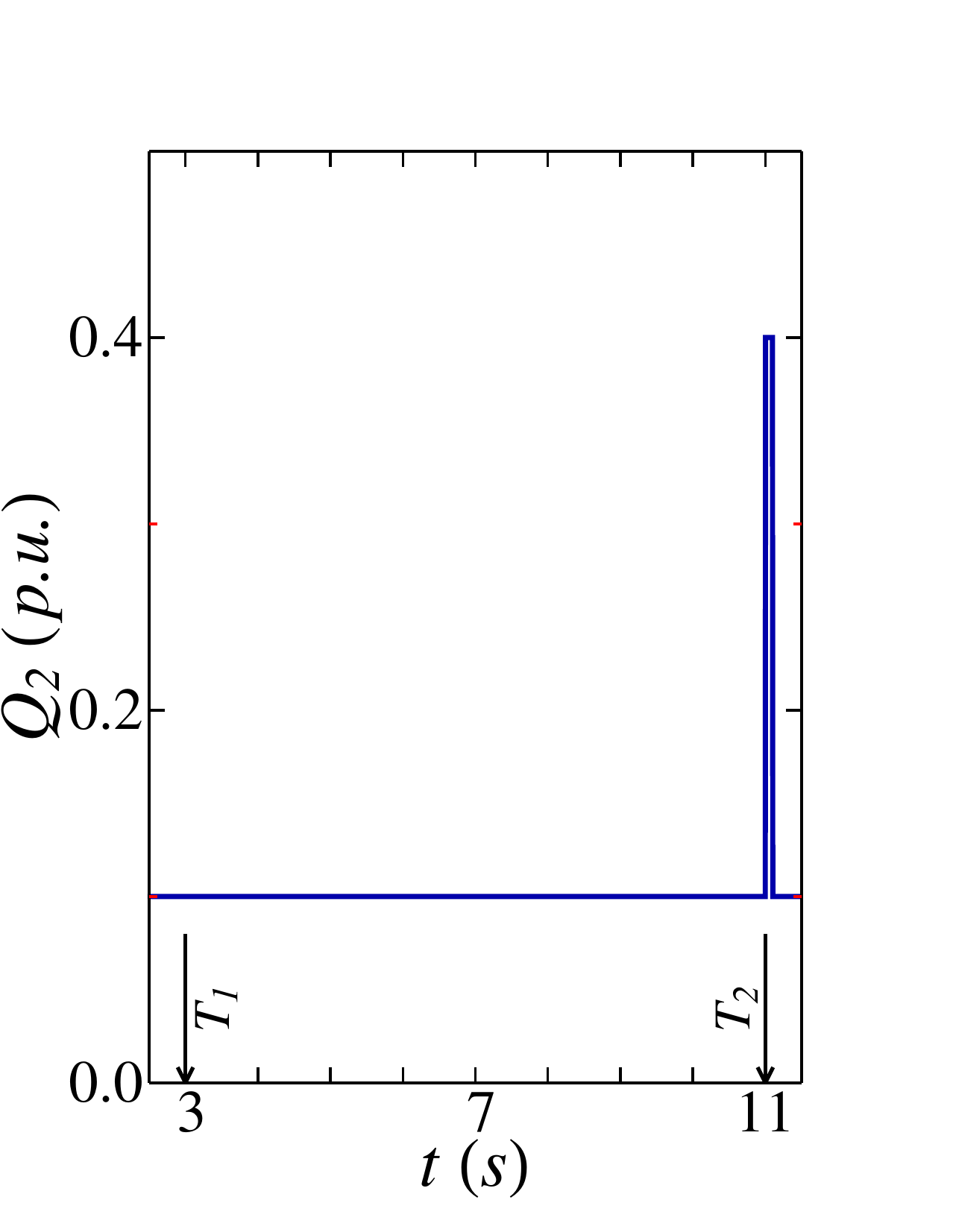} }}%
    \caption{Power demand level at load bus $2$, $t\leq{12} \,s$}%
    \label{fig:13busP2Q2}%
\end{figure}

\begin{figure}[!ht]
    \centering
    \includegraphics[width=0.8 \columnwidth]{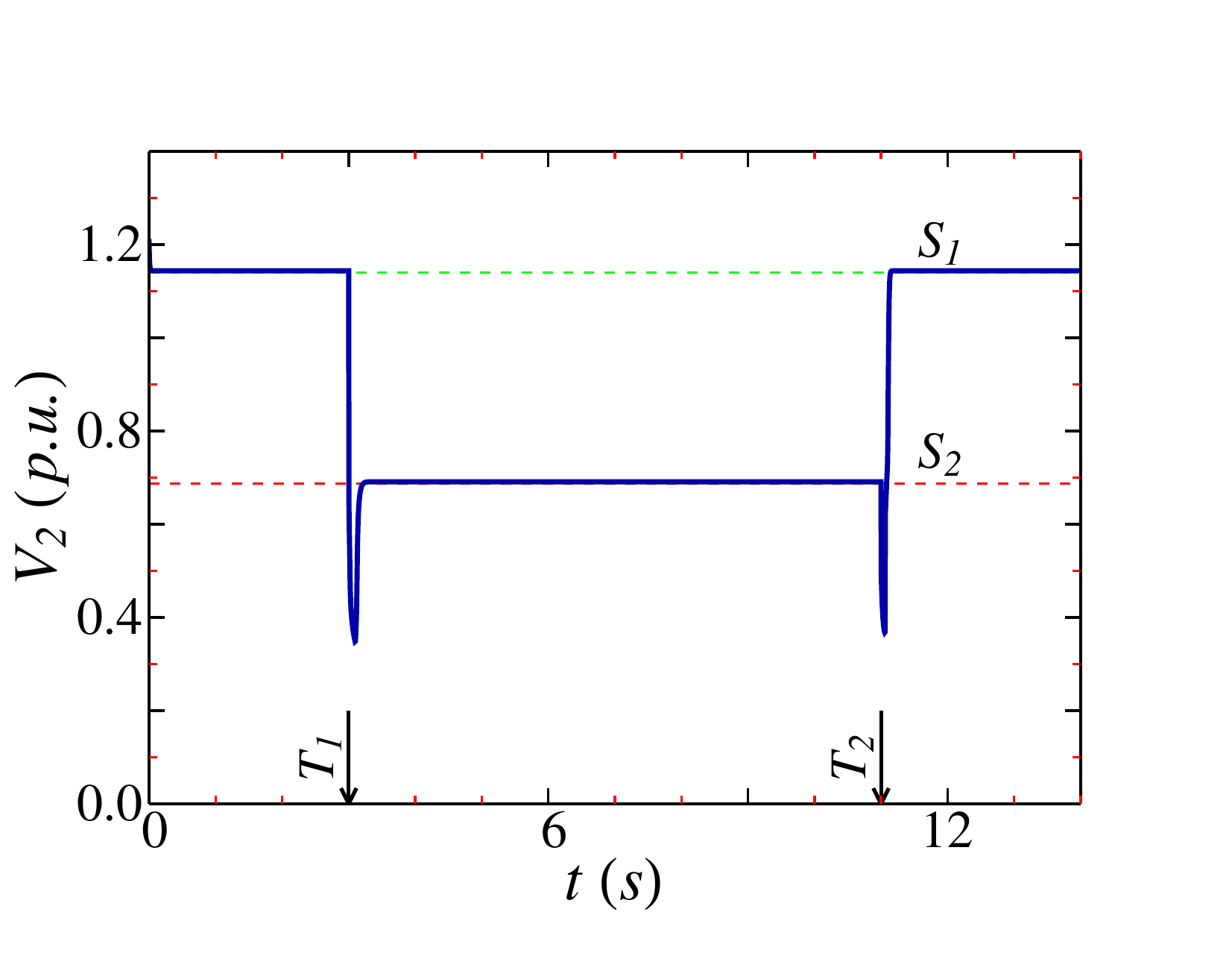}
	\caption{Voltage at load bus $2$}
    \label{fig:13busV2}
\end{figure}

To show that the system may get entrapped in the new equilibria, we conduct numerical experiment demonstrating the dynamic transitions between two stable equilibria in the $13$-bus system. We simulate the transient of the $13$-bus test case system between two stable equilibria $S_1$ and $S_2$ which are plotted in Figure \ref{fig:13busnose}. During the simulations, all loads except the load at bus $2$ are polynomial loads with fixed base demand levels, $P_{k\,0}$ and $Q_{k\,0}$. Load bus $2$ is composed of distributed generator and some dynamic load with aggregate time constants $\tau_1=\tau_2=0.01\,s$ and base demand levels, $P_2^s=-0.85\,p.u.$ and $Q_2^s=0.1\,p.u.$ Initially, the system is working at the high voltage equilibrium, $S_1$ and starts at the initial condition $g_2=0.077$, $b_2=-0.652$, $g_k=b_k=0$ for $k=\overline{3,13}$.

When the system is operating at the point $S_1$, a large disturbance occurs at $T_1=3\,s$ and lasts $0.1\,s$, during which the active power generation on bus $2$ is temporary lost, so the bus $2$ starts to consume active power with $P_2=0.5\,p.u.$ The disturbance is recorded in Figure \ref{fig:13busP2Q2}. After short period of time the generation on bus $2$ is automatically restored (after the action of re-closer or other automatic relay systems), and the power level returns back to the original levels. However, the system fails to restore to the normal operating point and becomes entrapped in the low voltage equilibrium. The restoration action that returns the system back to its normal operation conditions happens at $T_2=11 \,s$ and is described in details in next sections. To conclude this analysis, the loss of DGs results in a transient of the system from the high voltage equilibrium, $S_1$, to the lower voltage equilibrium, $S_2$, as shown in Figure \ref{fig:13busV2}.

\section{Discussion of the simulation results and the multistability phenomenon}

For the $3$-bus and $13$-bus systems considered in this simulation, the voltage of the second stable equilibrium is relatively low and unacceptable for the system to operate for a long period. However, even though the new solutions have unacceptably low voltage level for normal operation, they may have a strong effect on the process of post-fault voltage recovery where the voltage rises from low to normal values. This phenomenon is the main focus of this work and the main motivation to study the low voltage solutions and their stability. In other words, for voltage stability analysis purpose, we are interested in emergency situation rather than normal operation conditions.

In our numerical experiments, we did not observe any new solutions with acceptable levels of voltage. However, it is not clear whether this observation will hold in general for other types of distribution grids. For example, for a two-bus network with Under-Load Tap Changers (ULTCs) having limited taps, the voltage levels of both of the stable equilibria may be high enough for normal system operation which is considered in appendix \ref{app:ULTC}. Moreover, many of the current distribution grids do not have proper under-voltage protection on low voltage part. As long as the currents experienced during the nonlinear transients do not trigger the overcurrent relays the system may get entrapped at the low voltage equilibrium. In substations, system operators may not be aware that the system is working in such unfavorable conditions. Therefore, if no countermeasures are introduced, the system will stay at the low voltage equilibrium. 

Moreover, even for transmission grids, the system may experience post-fault low voltage for several seconds due to FIDVR. The low voltage conditions can move the system closer to the low voltage equilibrium and increase the probability of its entrapment. This may trigger the under-voltage protection relays and result in loss of service, and even cascading failures in most dramatic scenarios. In other words, the entrapment at the low voltage equilibrium may cause the same effects as FIDVR or power swings, i.e. the reliability risk of pro-longed low voltage conditions. To compare these phenomena, one can observe the similarities among post-fault voltage behaviors in Figure \ref{fig:shortcircuit} and Figure $1$, $3$, and $4$ in \cite{FIDVRwhitepaper}. However, multistability without countermeasures not only delays the voltage recovery process but also entraps the system at the low voltage equilibrium permanently. Therefore, multistabilty and the transition among equilibria need to study more seriously as FIDVR.

In future distribution networks, more electronic devices with wider operation ranges may be used to control the voltage of the loads. In this case, the voltage level on the consumer side may be in an acceptable range even when the voltage on the grid side is low.  At the same time, more PVs will be installed to supply power to individual consumers. If the ratio of distributed generation capacity to the total load capacity in the grid is large enough, which is may be as high as $35\%$ as recommended \cite{Wangthesis}, without the presence of energy storage, the randomness of weather conditions can cause large disturbance which result in the transient leads to system to be attracted at the low voltage equilibrium. To prevent the undesirable conditions, new monitoring and undesirable state detections schemes need to be introduced as well as additional preventive and corrective controls to keep the system in the normal operating point.

What happens if there is no multistability? By tuning the loads time constant in the simulations with the three-bus test system so that the low voltage equilibrium, $E_H$, is no longer stable, we re-run the simulation with the large disturbance at $t_{d1}$ to monitor the system behaviors. As expected, without the stable low voltage equilibrium, the system is simply recovered to the normal operation condition after being subjected to a severe disturbance. Therefore, multistablity prevents the normal restoration of the system after disturbances by entrapping the system at the low voltage equilibrium that also drives the system into emergency state. Proper emergency control actions need to restore the system to the normal operating conditions. For the restoration purpose, we design a new emergency control strategy in the next section.

\section{Pulse emergency control strategy-PECS} \label{sec:emergencycontrol}
\subsection{Existing emergency control schemes}

For small perturbations in distribution grid at low voltage which cause small deviations around the normal operating point, the $HV$ slack bus can strongly influence the voltage recovery process, thus reducing the need for load shedding. The support from $HV$ slack bus is provided first by instantaneously feeding of the low voltage grid with active and reactive power, and then, on longer timescales by adjustment of tap changers or voltage regulators, which is limited, say to $5\%$. In our simulations we assumed, that $HV$ slack bus may increase its voltage up to $1.05\, p.u.$; however, we observed that this effect is insufficient to recover voltage if the system is upset by large enough disturbances like the ones considered in the manuscript. This is so because the low voltage condition is not a result of the lacking reactive power support as in normal operating condition, rather the system is entrapped at another branch of the solution manifold following the transient response of the system. Therefore, injecting more reactive power into distribution system via Point of Common Coupling (PCC) does not necessarily restore the normal operating conditions. The same logic applies to the strategy of shedding the inductive loads and switching on capacitors banks. The existing voltage recovery systems are not guaranteed to work for non-regular operating points, where for example reduced reactive power consumptions does not necessarily lead to higher voltage level \cite{Overbye1994}. Hence, as shown in the transient among equilibria above, temporary loss of DG causes severe disturbances on the system and may force the entrapment of the system in the low voltage equilibrium; hence, causing the system to enter an emergency state.

In order to restore the system to the alert state, emergency control actions should be initiated \cite{Kundur}. For the distribution systems, appropriate emergency control actions may include fast capacitor switching, and fast load shedding. The performance of these control actions are assessed in this subsection. 

\subsubsection{Load shedding induced voltage collapse, LSIVC}
As the system enters emergency state, load shedding may be used in the last resort and seems to be one of the most reasonable and effective countermeasures. However in situations with reversed power flow, inappropriate load shedding strategies may result in voltage collapse, thus resulting in a novel phenomena named load shedding induced voltage collapse (LSIVC). This phenomenon may occur when the low voltage equilibrium is  close to the tip of the corresponding nose curve branch, or in other words, when the margin to the maximum generation point is very small. In this situations shedding of the load increases the generation levels on the system and forces the system beyond the maximum generation point, thus initiating the voltage collapse. 

\begin{figure}[!ht]
    \centering
    \includegraphics[width=1 \columnwidth]{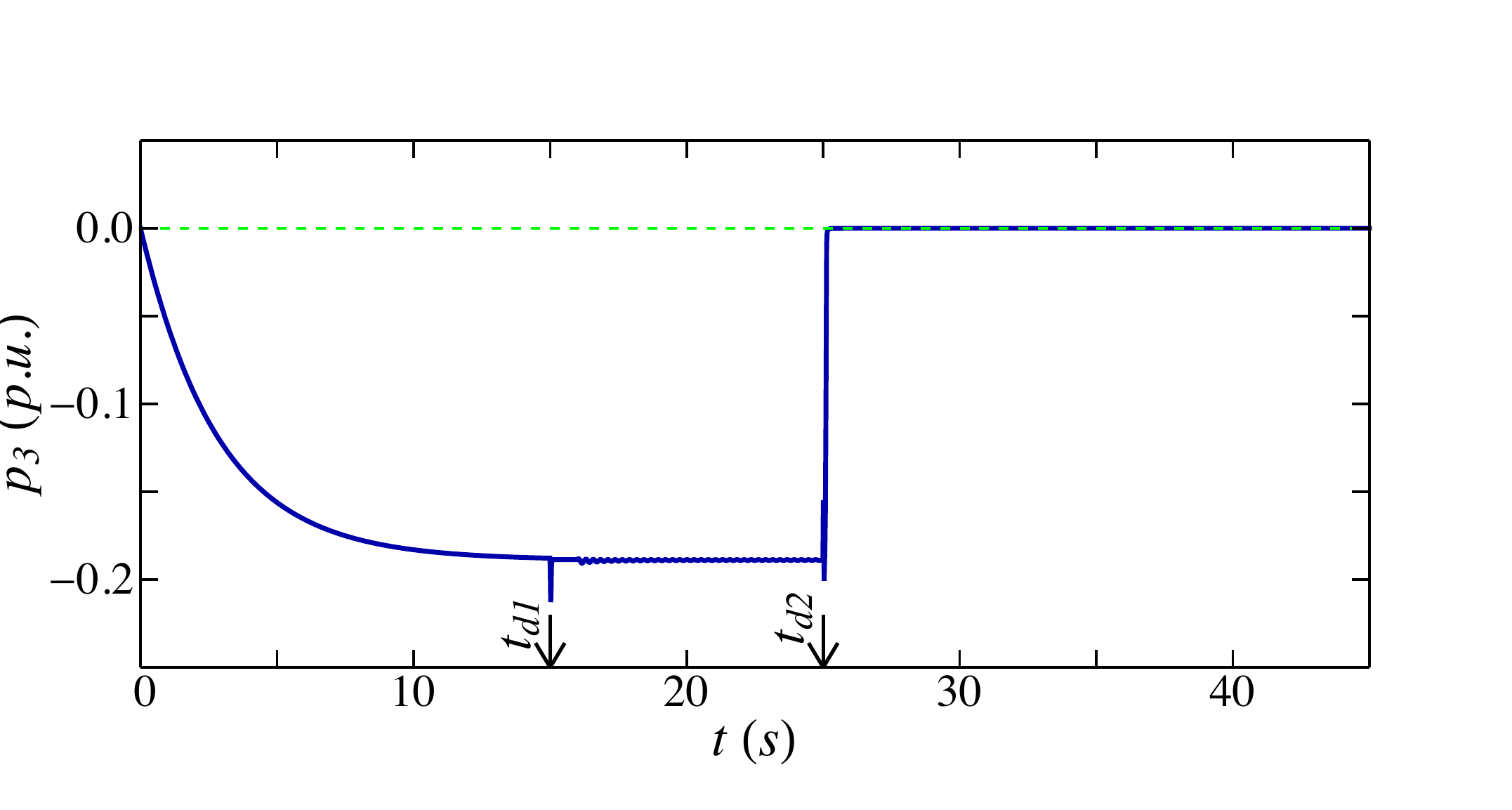}
	\caption{The active power of load bus $3$ in the load shedding induced voltage collapse}
    \label{fig:p3H2Lshed}
\end{figure}

For example, in the $3$-bus system described in a previous section, we consider a scenario where the system get entrapped in low voltage equilibrium with $P_2=- 0.7\, p.u$ after a transient fault shown in Figure \ref{fig:shortcircuit}. For this branch the tip of the nose corresponding to the maximal generation is at $P_2= - 0.73 \,p.u.$. In a scenario where the undervoltage protection system initiates shedding of the load at bus $2$, the overall generation at this bus is increased and the system goes beyond the maximal generation point. In a specific scenario, we have simulated shedding of $0.2 \,p.u$ at $t_{d2}=25\,s$ which resulted in a collapse of load $3$ so that the active power level remained zero after shedding load as shown in Figure \ref{fig:p3H2Lshed}.

Therefore, in the presence of multistability, load shedding cannot be always successfully employed for restoration of the system and undervoltage relays that initiate load shedding should be reconfigured. Same argument holds for the action of reactive power compensators that attempts to increase the load voltages by supporting more reactive power. The system may collapse as soon as it passes the tip of the $QV$ curves in the production region.

\subsection{Proposed pulse emergency control strategy}
To increase the reliability of the system without compromising the quality of power delivery service, we propose a new emergency control strategy that relies on DG curtailments. We observed that the system may converge to the undesirable low voltage equilibrium due to some pulse disturbance, like temporary loss of generation. Similarly, another pulse control force might help to recover the system to the normal operating condition. To force the system back to its normal operating point, the control action needs to temporarily move the system back to the regime without multistability, where the system can naturally converge to the normal operating point and naturally follow this branch returning back to the normal point when the pulse dies out. Thus, the problem is reduced to designing the appropriate duration, amplitude and composition of the pulse control action. Instead of providing precise instructions which requires rigorous transient dynamics analysis of the system, we propose the pulse choice heuristics based on the following logic. The magnitude of the pulse should be chosen in a way to drive the system out of the stability region of the undesired equilibrium point. Therefore, the magnitude of the pulse depends on the size of stability region of the low voltage equilibrium. An appropriate duration of the pulse should be long enough to allow the system to enter the stability region of the normal operating point. In a typical situation, the duration of the pulse should be compatible with the time constant of the dynamic loads, i.e. reasonably larger than the load time constant to provide enough time for the load to approach the new equilibrium. Our experiments also show that the faster the pulse is, the larger pulse magnitude is required. In addition, extremely large and long pulses may cause voltage collapse when the system moves outside the stability regions of both stable equilibria.

Normally, the distribution system without generation buses or PV buses has two solutions in the consumption regime where the loads consume power. For example, in Figure \ref{fig:voltagemultistability} and \ref{fig:13busnose}, the consumption regime is the right half plane where $P>0$. In most of the cases, the new solutions appear when the load starts exporting power into the network. However, it is possible for the system to have second solution branch even in the consumption regime. On the other hand, for normal grids, this branch can be observed in power consumption region only in relatively small neighborhood of zero consumption point. More thorough discussion is provided in the appendix \ref{app:newsol}. Therefore, if the control action could move the system to the consumption regime where only two solutions exist and mostly the high voltage solution is the unique stable one, the system will be driven to the upper branches of the nose curve which is characterized with high voltage level. After the pulse control action is cleared, the system will follow the upper branch to return to the high voltage equilibrium in the power generation regime. Based on this analysis, the proposed emergency control strategies dedicated to the DGs integrated distribution grids consist of two phases: the detection phase and the control phase. 

The objective of the detection phase is to identify the entrapment of the system in the undesirable low voltage equilibrium. This could be accomplished with the support of fast data acquisition systems such as Wide Area Measurement System/SCADA or future smart metering technologies. We envision a system where a large disturbance triggers the detector which then starts to monitor the system dynamics in an attempt to detect the events where the system enters in an emergency state where operating constraints are violated but the load constraints are still satisfied. As soon as the transition between equilibria is confirmed, the control phase is initiated. Large capacity DGs are chosen as candidates for PECS. A suggested set of candidates includes those have initiated the  initial disturbance leading to the entrapment. The DGs candidates are then curtailed for a short period of time, comparable to the natural relaxation times of the system. The DGs curtailments will drive the system to the consumption regime where the system recovers to the normal operating branch.

In some situations, the pulse emergency control design based on DGs curtailments may need more precise and rigorous calculation and analysis. Hence, we provide an alternative design of emergency control system. In this alternative the control directly changes the set of impedance states, i.e. the conductances $g$ and susceptances $b$ of the loads, rather than load powers or voltage levels. This allows the controller to move the system state in the stability region, from which it will provably reach the desired equilibrium eventually. In order to illustrate this control technique, we modify slightly the dynamic load models as below.

\begin{align}
\tau_g\,\dot{g}=-(g-g^s)\\
\tau_b\,\dot{b}=-(b-b^s)
\end{align}

where $g^s$ and $b^s$ are the set points of impedance which may be set to be equal the state of the desired equilibrium which is the normal operating point in this case; $\tau_g$ and $\tau_b$ are the time constants compatible to the load time constants in \eqref{eq:gg} and \eqref{eq:bb}. Note, that the modified form of the load dynamics is equivalent to the one described in \eqref{eq:gg} and \eqref{eq:bb} and the only difference is that the modified form applies to the control period. By using this type of controller and appropriate set points $g^s$ and $b^s$, a long enough PECS period will allow the system enter the stability region of the high voltage equilibrium. Then the system successfully returns to the normal operating point.

The phase portrait also helps to understand the trajectory of the system and the transients among equilibria. It may also indicate which disturbances cause the system to ``jump''.

\begin{figure}[ht]
    \centering
    \includegraphics[width=0.75 \columnwidth]{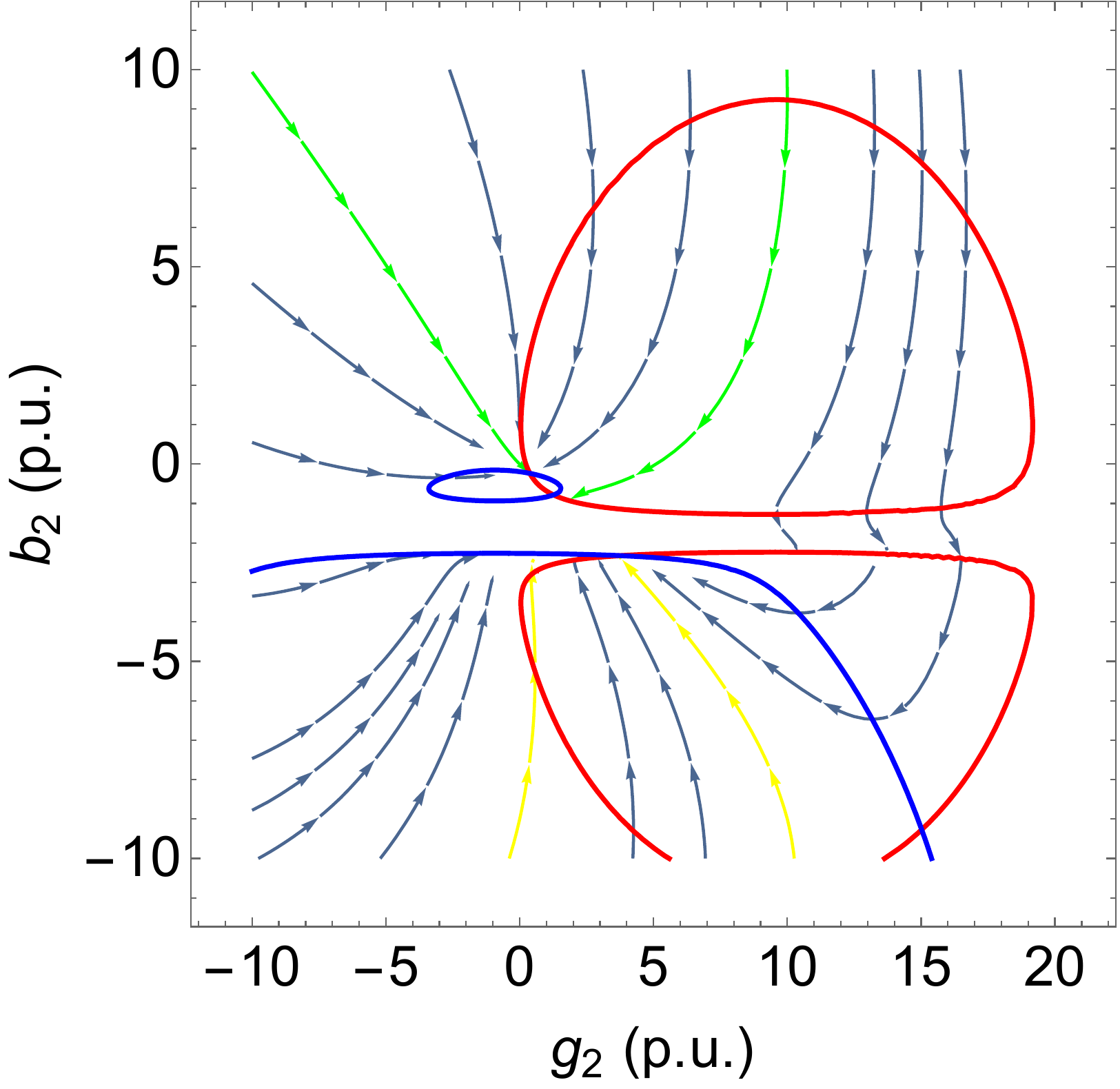}
	\caption{The phase portrait of the $3$-bus system}
    \label{fig:phaseportrait}
\end{figure}

The stability region which corresponds to a stable equilibrium can be visualized via phase portraits. For simplicity, we illustrate this with a $2D$- phase portrait for a simplified $3$-bus system as shown in Figure \ref{fig:phaseportrait}. The  test case parameters are as the following $z_{12} = z_{23} = 0.095+j 0.448$, $P_2 = 0.235$, $Q_2 = -0.145$, $g_3 = -0.246$, $b_3 = 1.46\, b_2$, $\tau_{1}= 0.56$, $\tau_{2}= 0.489$. All units are in $p.u.$. In Figure \ref{fig:phaseportrait}, $4$ solutions to the power flow problem are the intersection between red curves and blue curves. The state velocity vector field shows that there are two stable equilibria. One can observe that both stable equilibria have relatively large attraction regions. The attraction region of high voltage equilibrium is then used to design emergency control action that ensures that the system is driven into the basin of attraction.

Certainly, the phase portrait approach will only work for simple systems in practice, and more sophisticated algorithms are needed to characterize the attraction region. Construction of such algorithms is tightly linked to the problem of transient stability and is well beyond the scope of the present study. 

\begin{figure}[ht]
    \centering
    \includegraphics[width=0.7 \columnwidth]{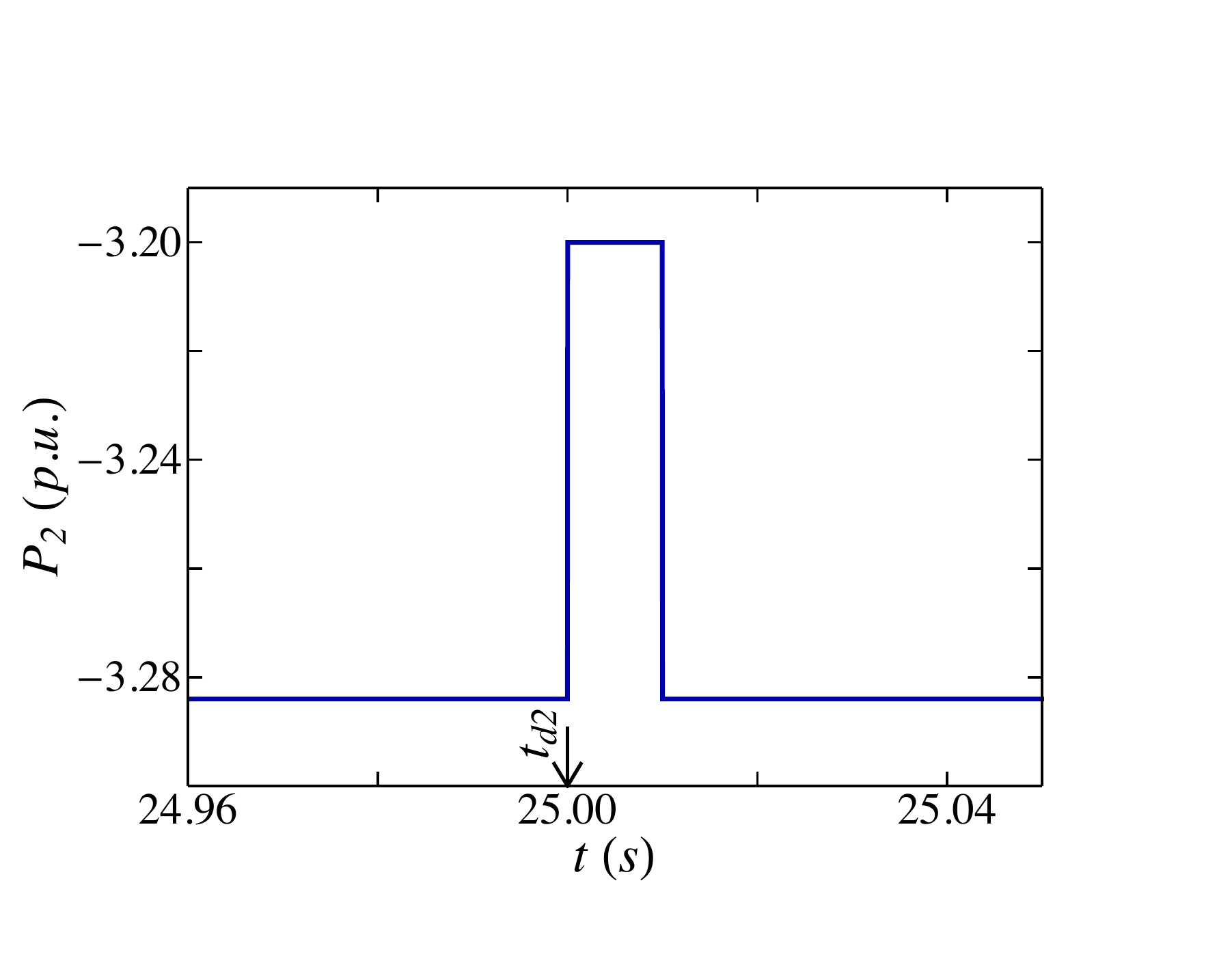}
	\caption{Active power at bus $2$ in PECS}
    \label{fig:L2HP}
\end{figure}

\begin{figure}[htb]
    \centering
    \subfigure[]{{\includegraphics[width=4.25cm]{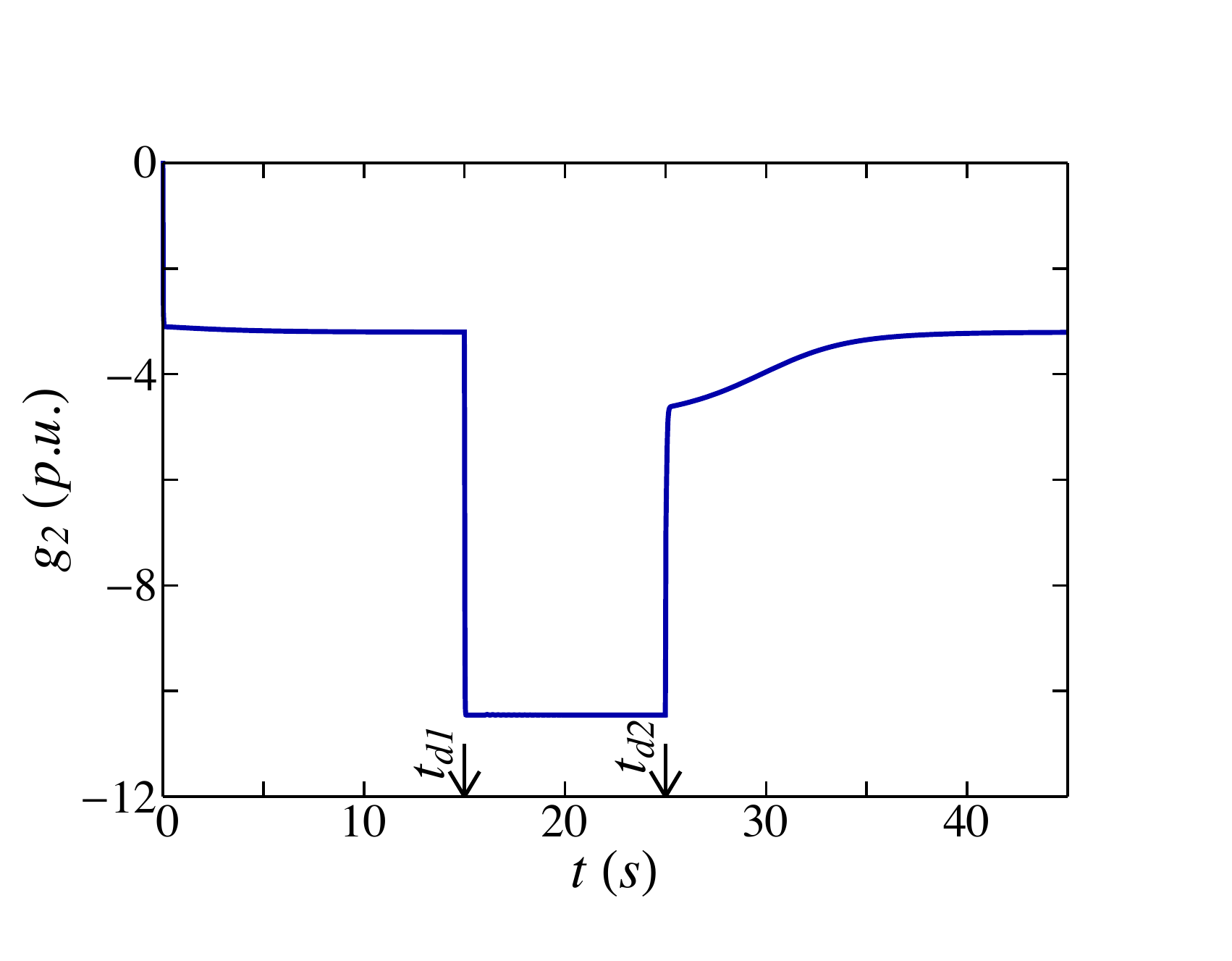}}}%
    \subfigure[]{{\includegraphics[width=4.25cm]{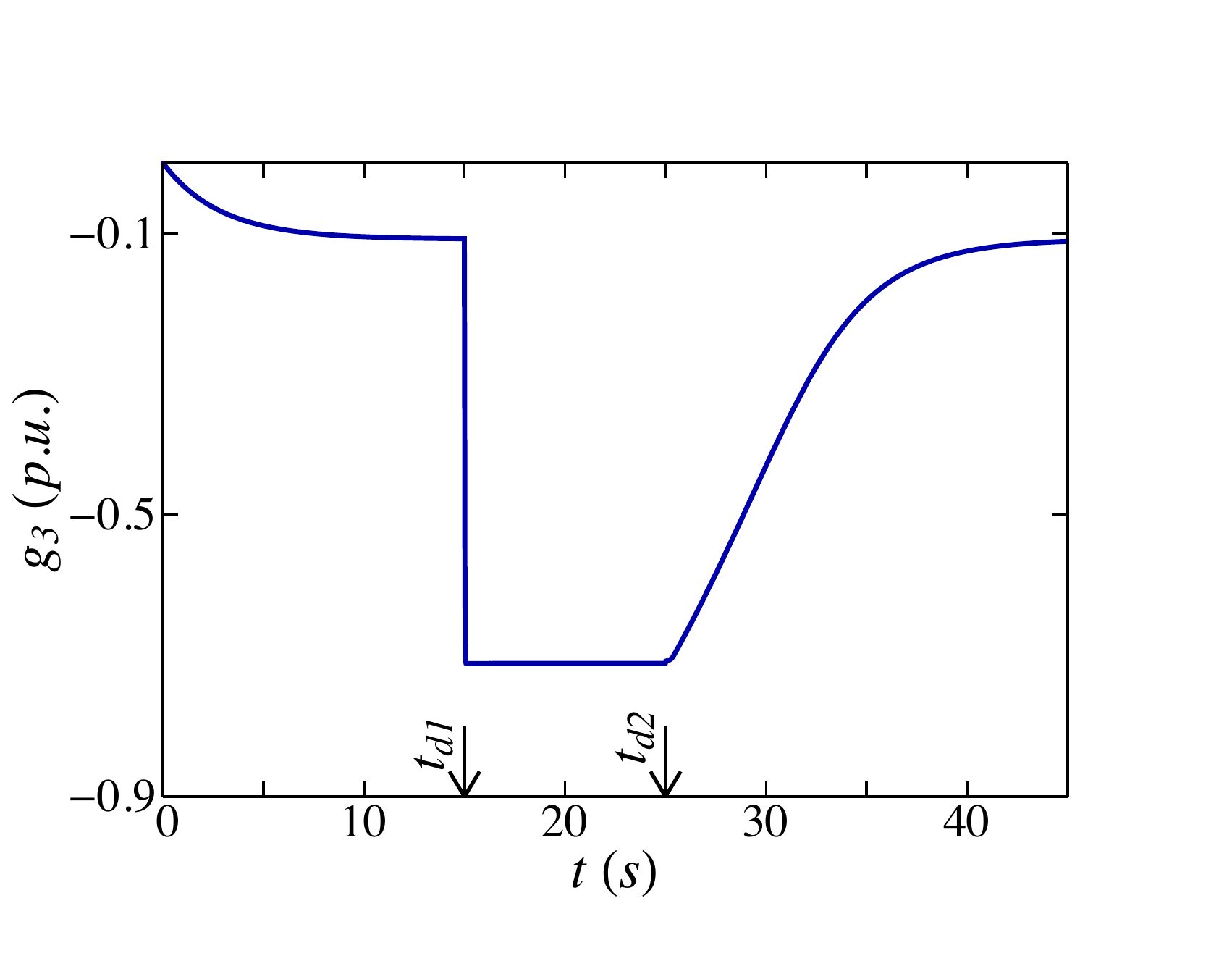} }}%
    \quad
    \subfigure[]{{\includegraphics[width=4.25cm]{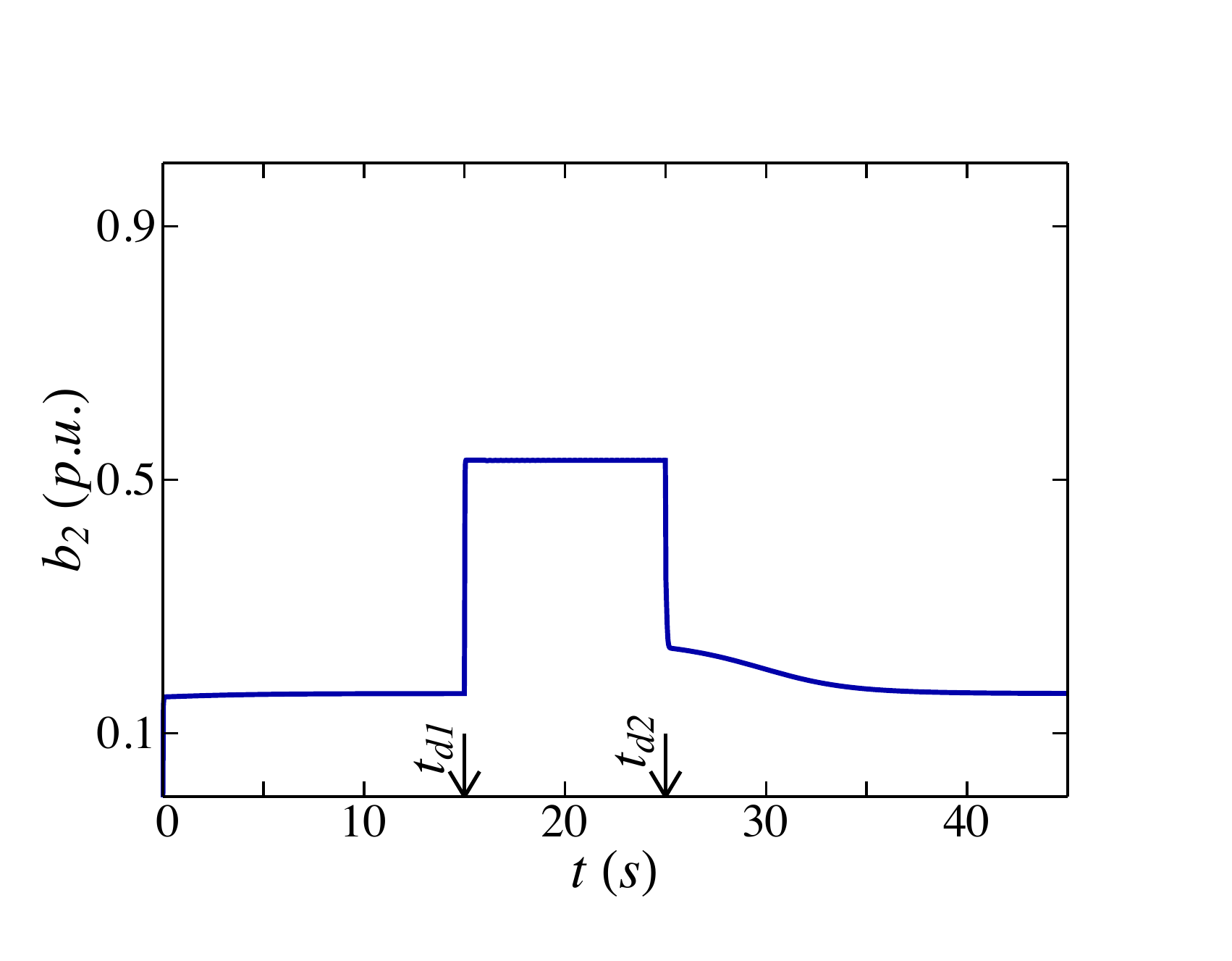}}}%
    \subfigure[]{{\includegraphics[width=4.25cm]{g3_3bus.pdf} }}%
    \caption{The conductances and susceptances at bus $ 2$ and bus $3$, $ t\leq{45} \,s$}%
    \label{fig:gbL2H}%
\end{figure}

\begin{figure}[ht]
    \centering
    \includegraphics[width=0.8 \columnwidth]{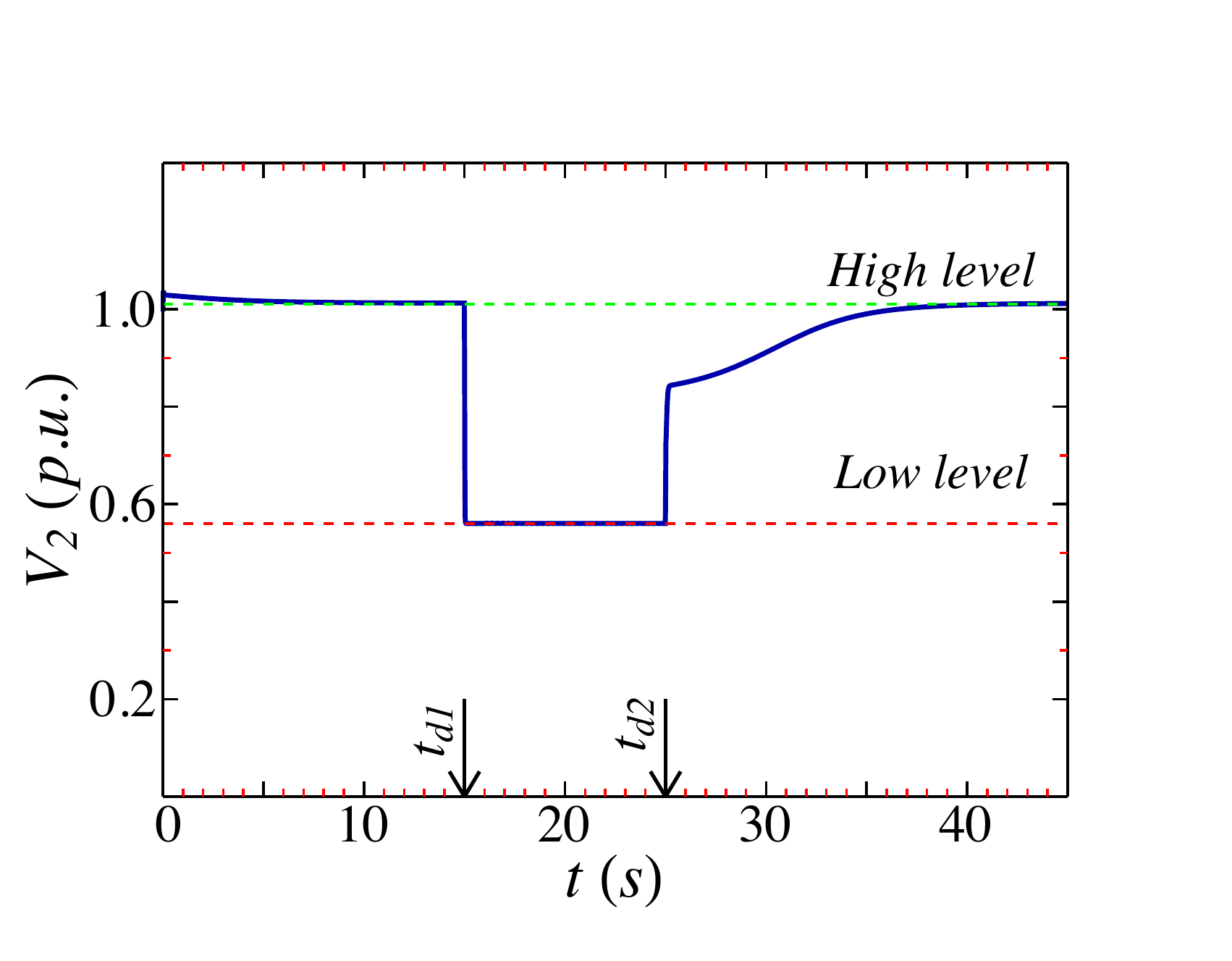}
	\caption{Voltage at bus $2$ for the large second disturbance, $t\leq{45} \,s$ }
    \label{fig:v2L2H}
\end{figure}

We implemented the proposed emergency control strategy both for the two test cases described in this work. For the $3$-bus system, the pulse emergency control action is applied to restore the system to the high voltage stable equilibrium, $E_H$, after being entrapped at the low voltage stable equilibrium, $E_L$. As was shown in the section \ref{H2L}, the system is entrapped at the low voltage equilibrium at $t < {15.3\,s}$. The pulse emergency control action occurs at $t=25\, s$ and causes active power demand level at bus $2$, $P_2^s$, to increase for $0.01\,s$ as depicted in Figure \ref{fig:L2HP}. Consequently, the system returns to the high voltage equilibrium with the trajectory shown in Figure \ref{fig:gbL2H} and Figure \ref{fig:v2L2H}.

For the $13$-bus feeder system, a pulse emergency control action is also taken place to restore the system from the low voltage equilibrium, $S_2$. The candidate for the curtailment is the load bus $2$ as shown in Figure \ref{fig:13busP2Q2} at $T_2=11\,s$. As a result, the system is successfully recovered. The restoration of the system is recorded in Figure \ref{fig:13busV2}.

Moreover, the operator should be notified whenever possible about the DGs step reduction in advance. The optimal countermeasures should be designed in cooperation between transmission and distribution grid control systems. As we observed in our test case, the impact of slack bus or PCC connected to transmission grid is limited, especially in transient response of the system, and is insufficient to recover system from the low voltage level of the second equilibrium. It's worth to emphasize that the low voltage condition when the system is entrapped at the undesirable equilibrium is not a result of lacking reactive power support, rather it is so due to the existence of unusual solution branches on the solution manifold. At the same time, the transient among equilibria may happen as a result of transient faults which cannot be anticipated. As a matter of fact, the observed low voltage levels are a result of post-fault transient recovery. Obviously, for unexpected faults, coordination and anticipation is not always possible.

To prevent the entrapment of the system at the undesirable states new policies for power reversal need to introduced. The stability of the system depends both on the active and reactive power dynamics, so the regulations should be based on the analysis of the accurate models of distribution system dynamics. Standardization based on power factor may not guarantee the stability of higher voltage branches. The existing standards for DG penetration may not be adequately assessing the voltage reliability \cite{Wangthesis} and security of the system. Unlike transmission grids, the distribution systems are usually operated without designated distribution system operator monitoring the state of the grid and usually rely on the fully automated control. This situation is unlikely to change in the nearest future, and thus having more advanced automatic control systems capable of detecting the entrapment of the system and designing the optimal emergency control actions is necessary for reliable operation of future power grids with high penetration of distributed renewable or gas-based generation.

\section{Conclusion}
In this work we have shown that distribution grids with active or reactive power flow reversal can have multiple stable equilibria. The transitions among equilibria are also possible and occur after common disturbances. We demonstrate the existence and stability of the new solutions by considering several scenarios on a simple three-bus and more realistic $13$-bus feeder models. To identify the new branches we introduce a novel admittance homotopy methods which has better convergence properties in comparison to more traditional alternatives. Existing emergency control actions may fail to restore the system back to the normal operating conditions, and may even aggravate the situation. We demonstrate this phenomenon with an example of Load Shedding Induced Voltage Collapse that may occur in the system exporting active power. To address this problem we propose a novel emergency control scheme called Pulse Emergency Control Strategy. This strategy can successfully restore the system with temporary curtailment of the distributed generation. However, more studies are needed to determine the optimal amount and time of control actions. 

\section{Acknowledgement}
The work was partially supported by NSF, MIT/Skoltech and Masdar initiatives, Vietnam Educational Foundation, and the Ministry of Education and Science of Russian Federation, Grant Agreement no. 14.615.21.0001.

\appendices

\section{Multistability in a ULTC system} \label{app:ULTC}
\begin{figure}[ht]
    \centering
    \includegraphics[width=0.8 \columnwidth]{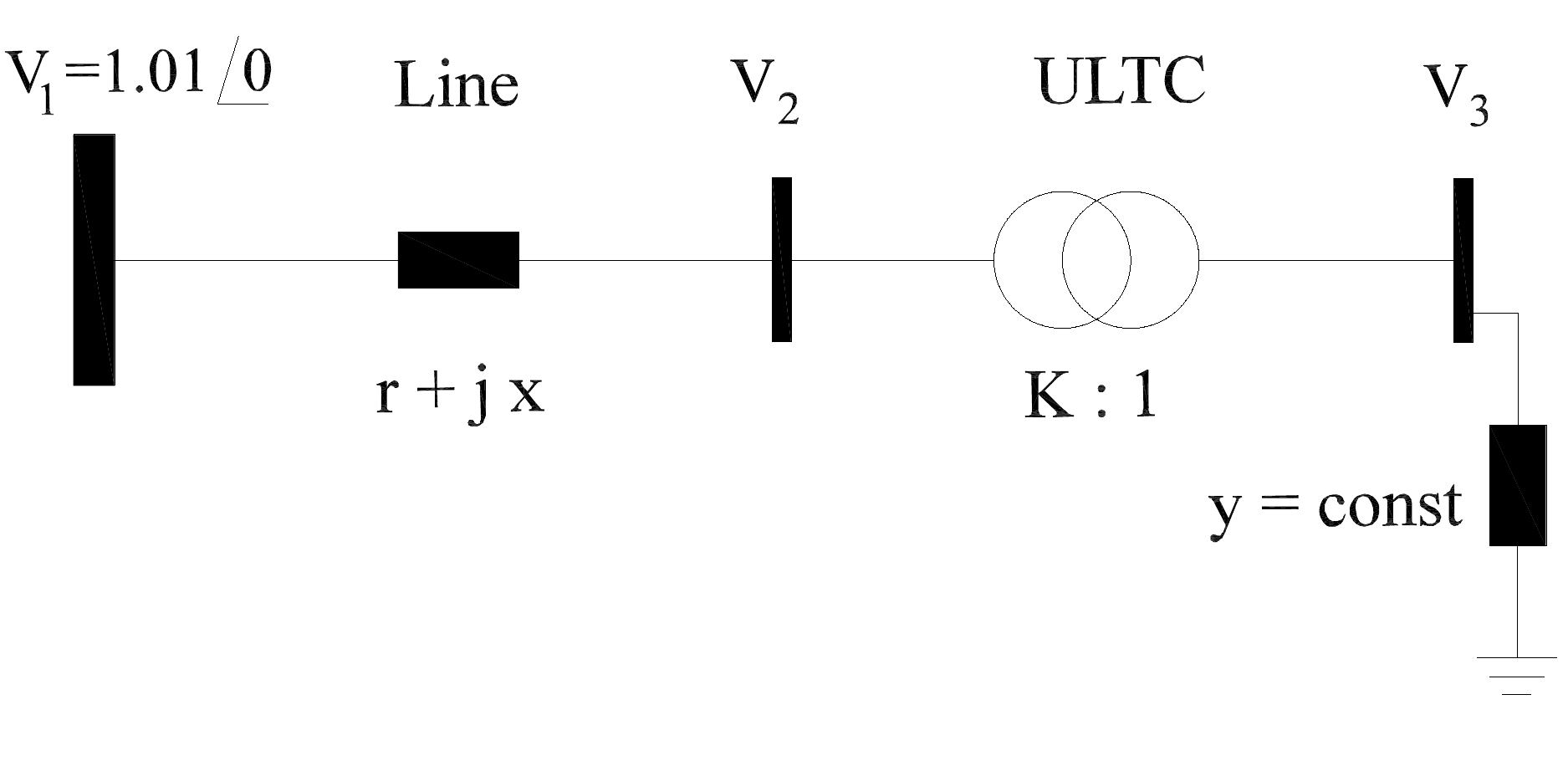}
	\caption{A $3$-bus system equipped with an ULTC}
    \label{fig:ULTCnet}
\end{figure}

In this section, we demonstrate how voltage-multistability can appear in normal operating conditions in a $3$-bus system where the voltage is controlled by the standard ULTC as depicted in Figure \ref{fig:ULTCnet}. We show that in this network there are two stable equilibria which both lie within the acceptable range of voltage.

\begin{figure}[ht]
    \centering
    \includegraphics[width=0.7 \columnwidth]{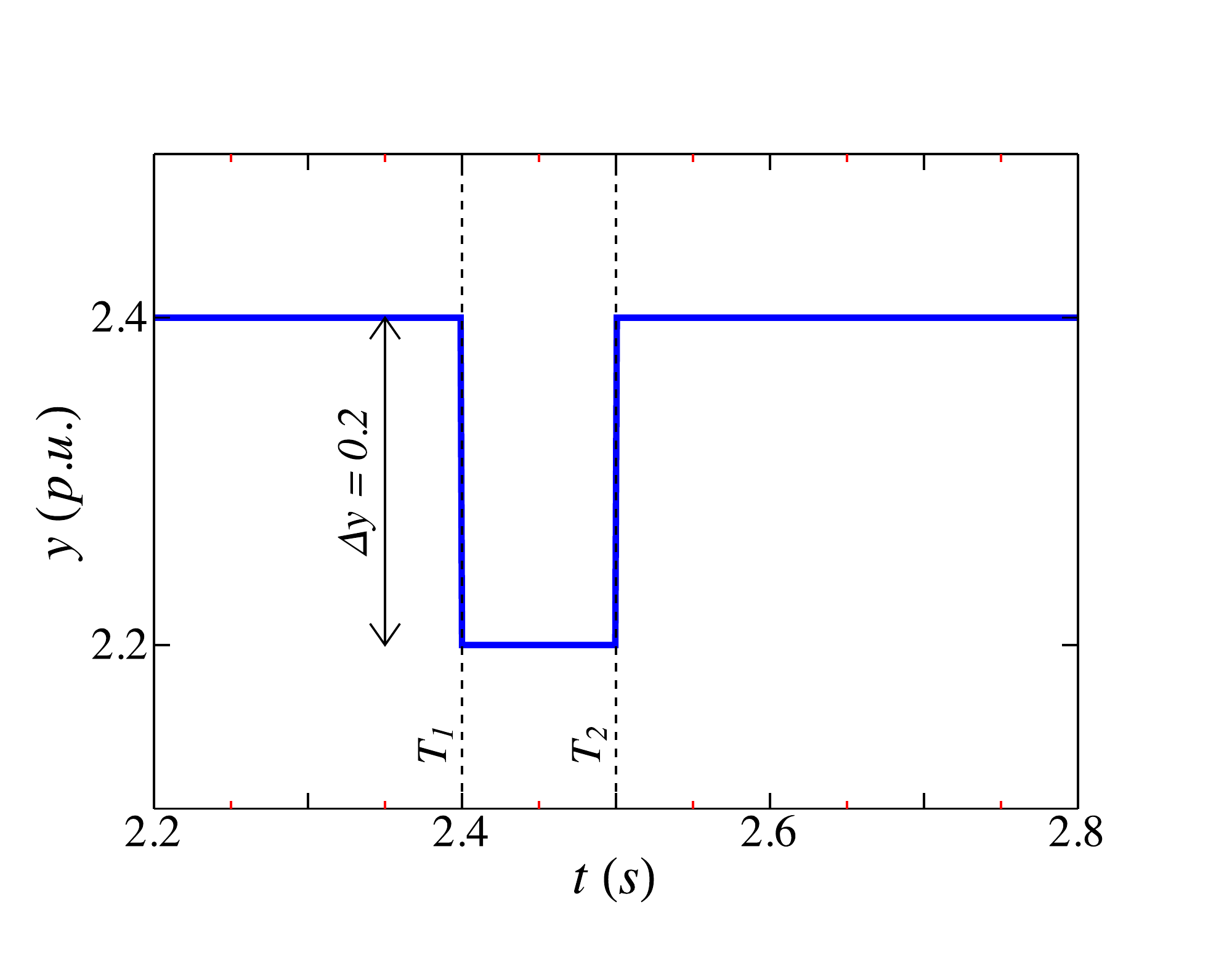}
	\caption{The admittance of the load during the disturbance}
    \label{fig:ULTCy}
\end{figure}

The network and the ULTC parameters are as follows:
ULTC deadband, $DB= 1.5\,\%$, i.e. $V^{max} = 1.015\, p.u.$, $V^{min} = 0.985\, p.u.$; $K^{max} = 1.17$, $K^{min} = 0.83$; $r = 0.069\, p.u.$; $x = 0.258 \,p.u.$; $cos\phi = 0.77$; $V_1 = 1.01\, p.u.$. The considered impulsive disturbance has the magnitude of $\Delta{y} = 0.2 \,p.u.$ and duration of $\Delta{t} = 0.1 \,s$ as shown in Figure \ref{fig:ULTCy}. Also, we use the continuous model to describe ULTC, hence ULTC control scheme can be modeled as:

\begin{equation}
\dot{K} = 
\begin{cases} 
-1, & \mbox{if } V_{3}<V^{min}\, \mbox{and}\, K>K^{min}\\ 
+1, & \mbox{if } V_{3}>V^{max}\, \mbox{and}\, K<K^{max} \\
\,0, & \mbox{otherwise}
\end{cases}
\end{equation}

The continuous time model here is chosen for the sake of simplicity, similar phenomena also can be observed with discrete model. More details of ULTC modeling can be found in \cite{milano2011hybrid}.

\begin{figure}[ht]
    \centering
    \includegraphics[width=0.8 \columnwidth]{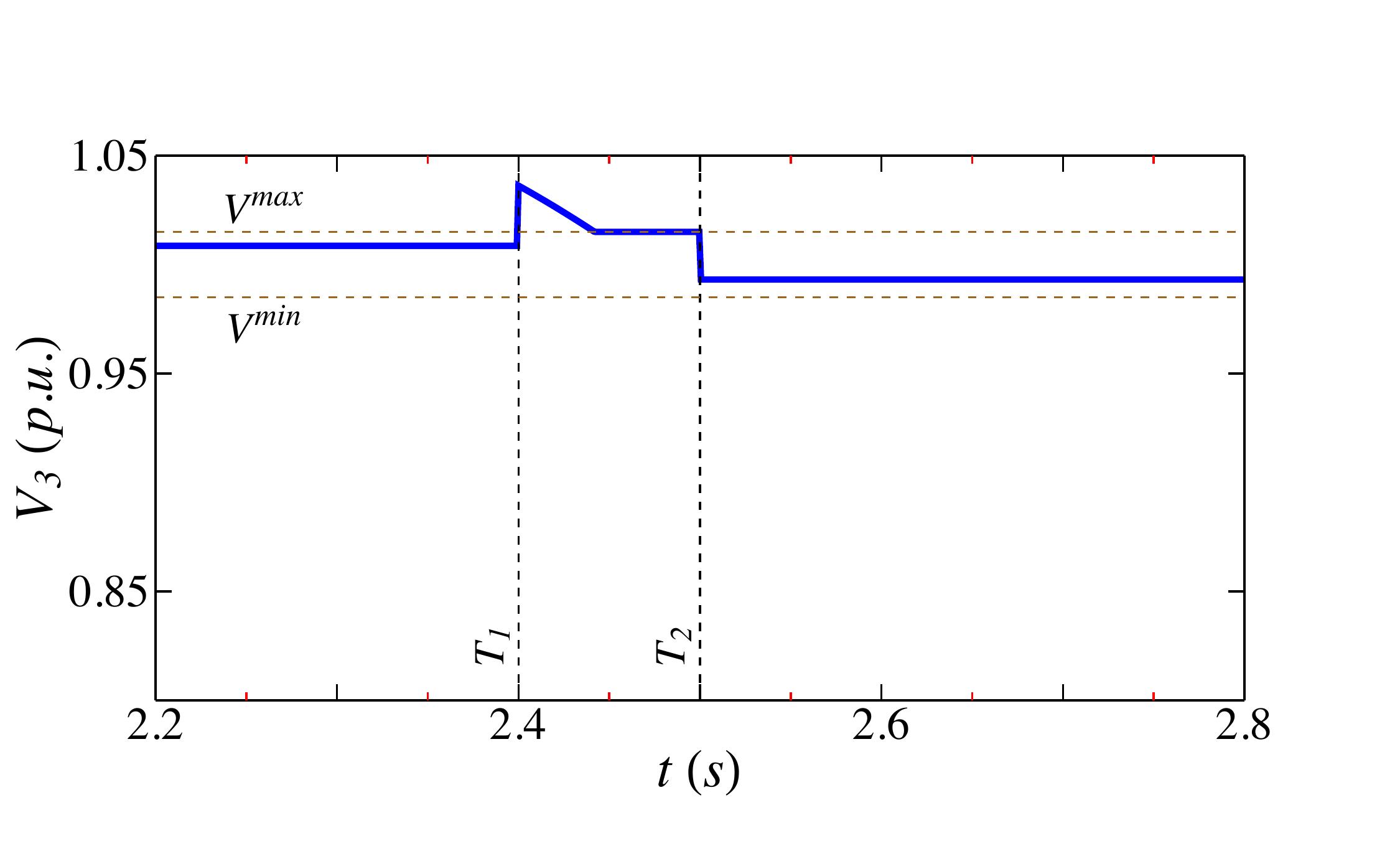}
	\caption{Voltage multistability with the ULTC system}
    \label{fig:ULTCV}
\end{figure}

Apparently, two different multi-stable operating points shown in Figure \ref{fig:ULTCV} which correspond to either $V_3 = 1.009 \, p.u.$ or $V_3 = 0.993 \,p.u.$, satisfy voltage constraints.

In the absence of ULTC the constant impedance load system is linear, and possesses only one equilibrium for a given impedance value. The ULTC plays the role of an effective load controller, but unlike standard dynamic models of $PQ$ type loads that have fixed power consumption level at equilibrium, the ULTC has an controls the level of voltage and can have multiple equilibria either when there is a deadband, or when the tap ratio hits the limits. This simple example is presented just to illustrate the possible sources of multistability phenomena in power systems. This example was illustrating how multistability appears in the presence of controller deadbands, below we also show that the new equilibria can appear in the presence of tap ratio limits.

\begin{figure}[ht]
    \centering
    \includegraphics[width=0.8 \columnwidth]{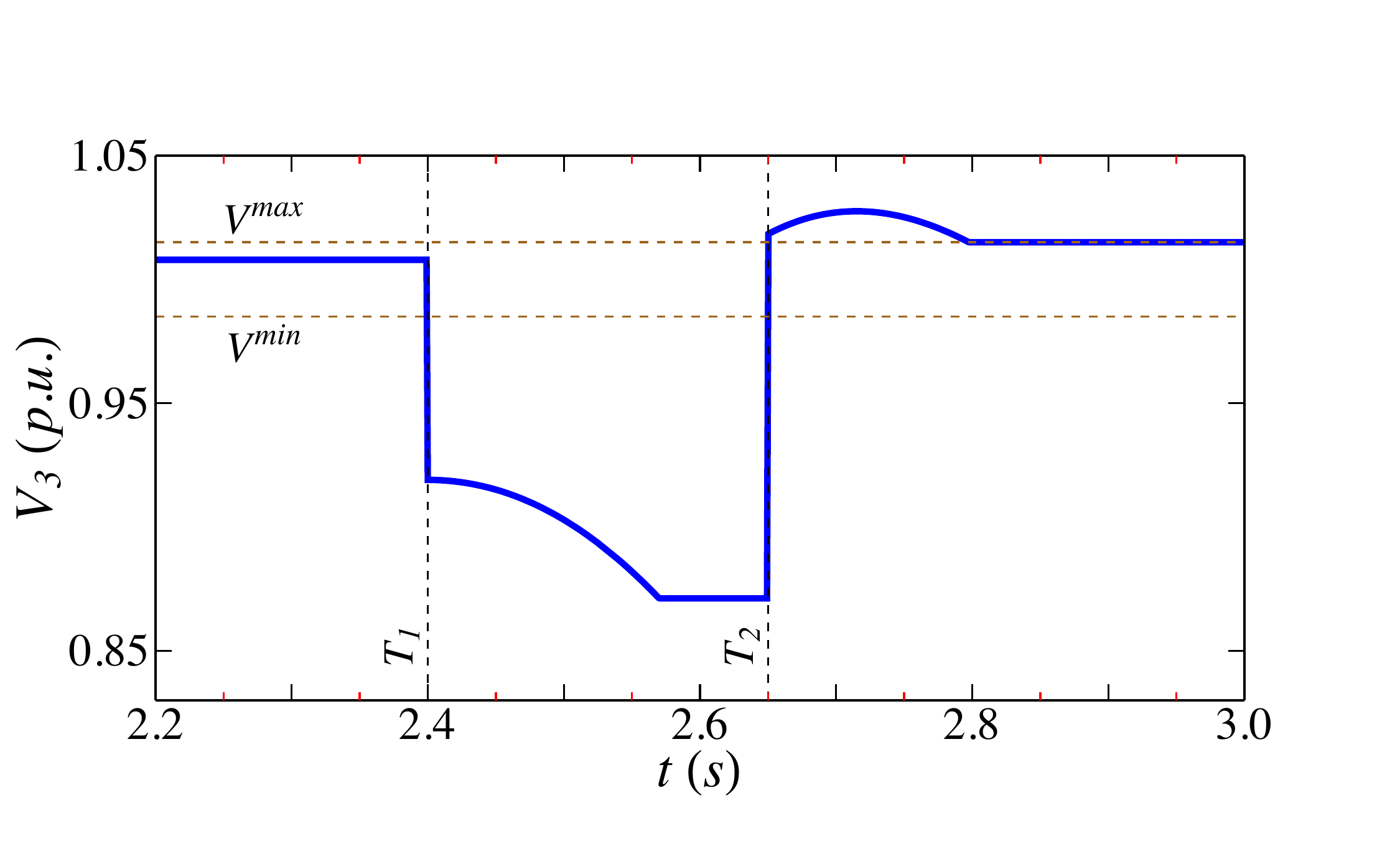}
	\caption{Voltage multistability with the ULTC system due to tap limits}
    \label{fig:ULTCV_hys}
\end{figure}

\begin{figure}[ht]
    \centering
    \includegraphics[width=0.8 \columnwidth]{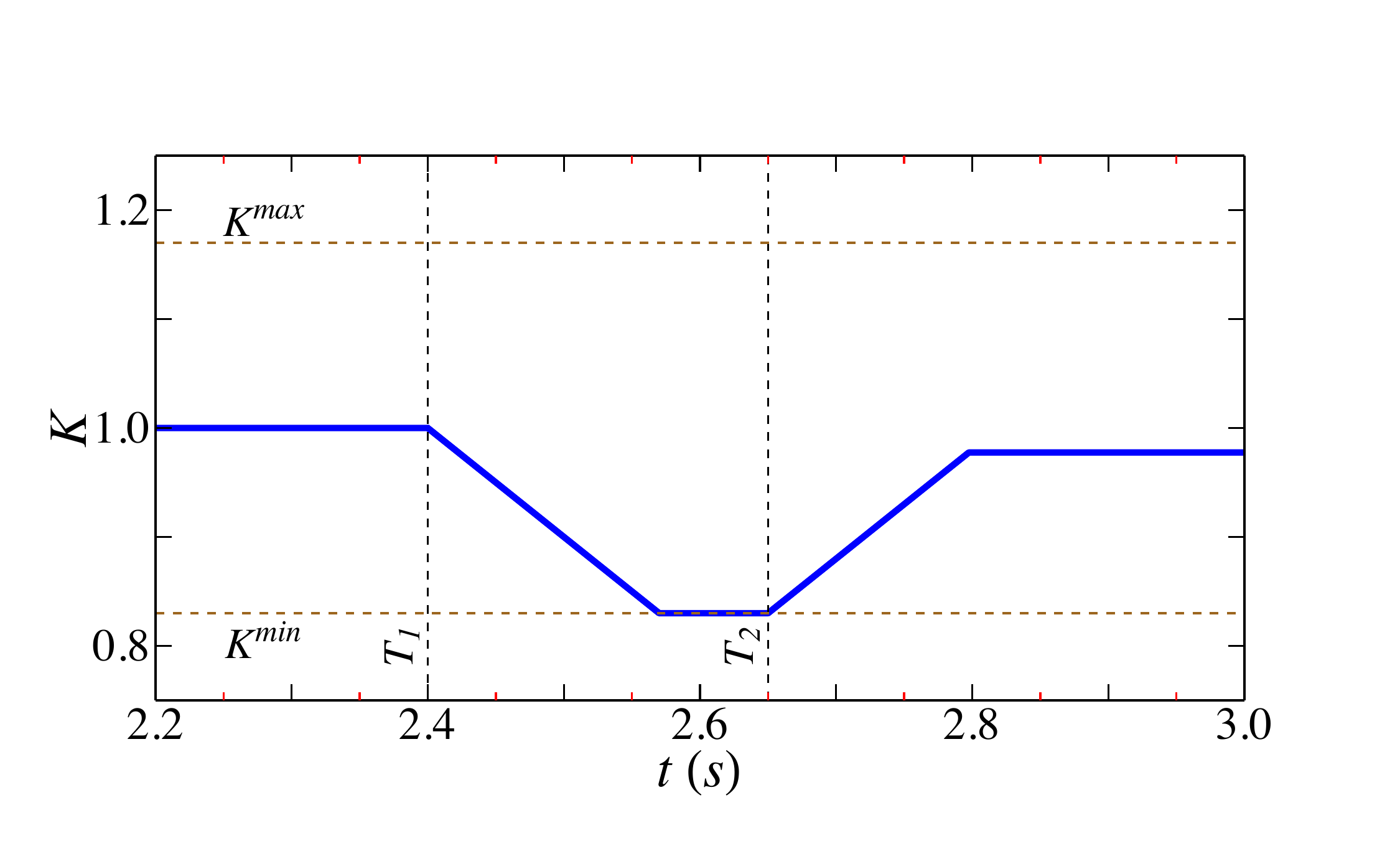}
	\caption{The ULTC ratio during the disturbance}
    \label{fig:ULTCy_hys}
\end{figure}

To illustrate the effect of tap limits, we reconsider the ULTC test case with new parameters are as follows: $\cos\phi = 0.8$; $V_1 = 1.03\, p.u.$. Other parameters are unchanged. The impulsive disturbance has the magnitude of $\Delta{y} = -0.6 \,p.u.$ and duration of $\Delta{t} = 0.25 \,s$ between $T_1 = 2.4\,s$ and $T_2=2.65\,s$. Multistability is then observed in either Figure \ref{fig:ULTCV_hys} or Figure \ref{fig:ULTCy_hys}. Figure \ref{fig:ULTCy_hys} shows that the tap ratio lower limit $K^{min}$ is reached. When the load impedance returns to the predisturbance level at $t=T_2$, the tap ratio converges to the post-disturbance value which differs from the predisturbance one. Therefore, the tap limits also can lead to multistability.

\section{Appearance of new solution branches} \label{app:newsol}

The appearance of new solutions in power flow reversal regime is not a mathematical fact, but rather an empirical observation, that can be justified by some qualitative reasoning. It is possible for the system to have second solution branch even in consumption regime; but for normal grids, this branch can be observed in power consumption region only in relatively small neighborhood of zero consumption point. We illustrate this claim by considering a simple $3$-bus test case as described below. Although this example is overly simplistic, the behavior observed is qualitatively similar to other more realistic systems. Despite extensive testing on several IEEE cases, we did not observe the non-trivial branch in power consumption regime for operating conditions corresponding to normal loading levels.

\begin{figure}[ht]
    \centering
    \includegraphics[width=0.8\columnwidth]{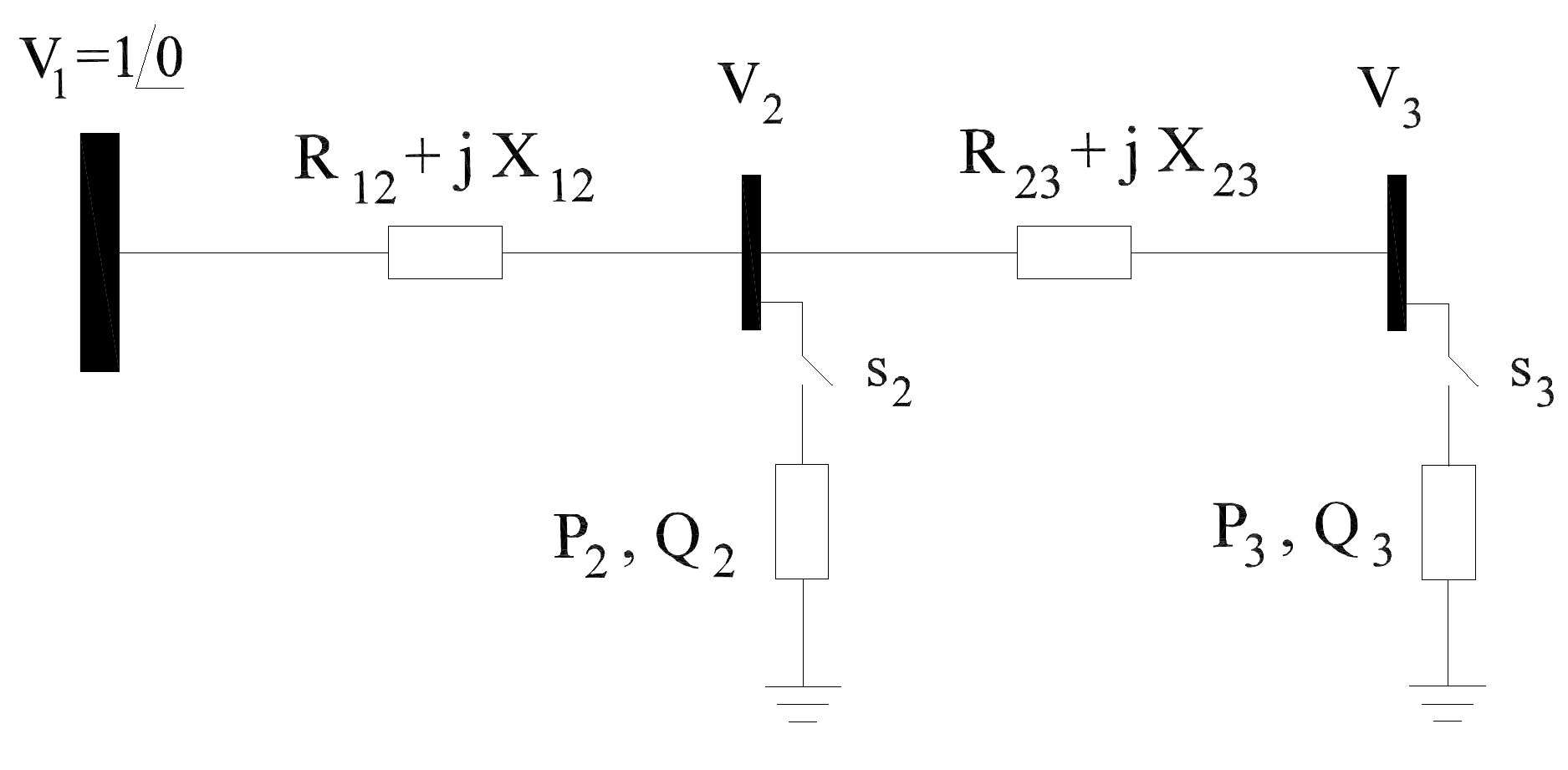}
	\caption{A $3$-bus system with fictitious switches at load buses}
    \label{fig:3busSW}
\end{figure}

\begin{figure}[ht]
    \centering
    \includegraphics[width=1\columnwidth]{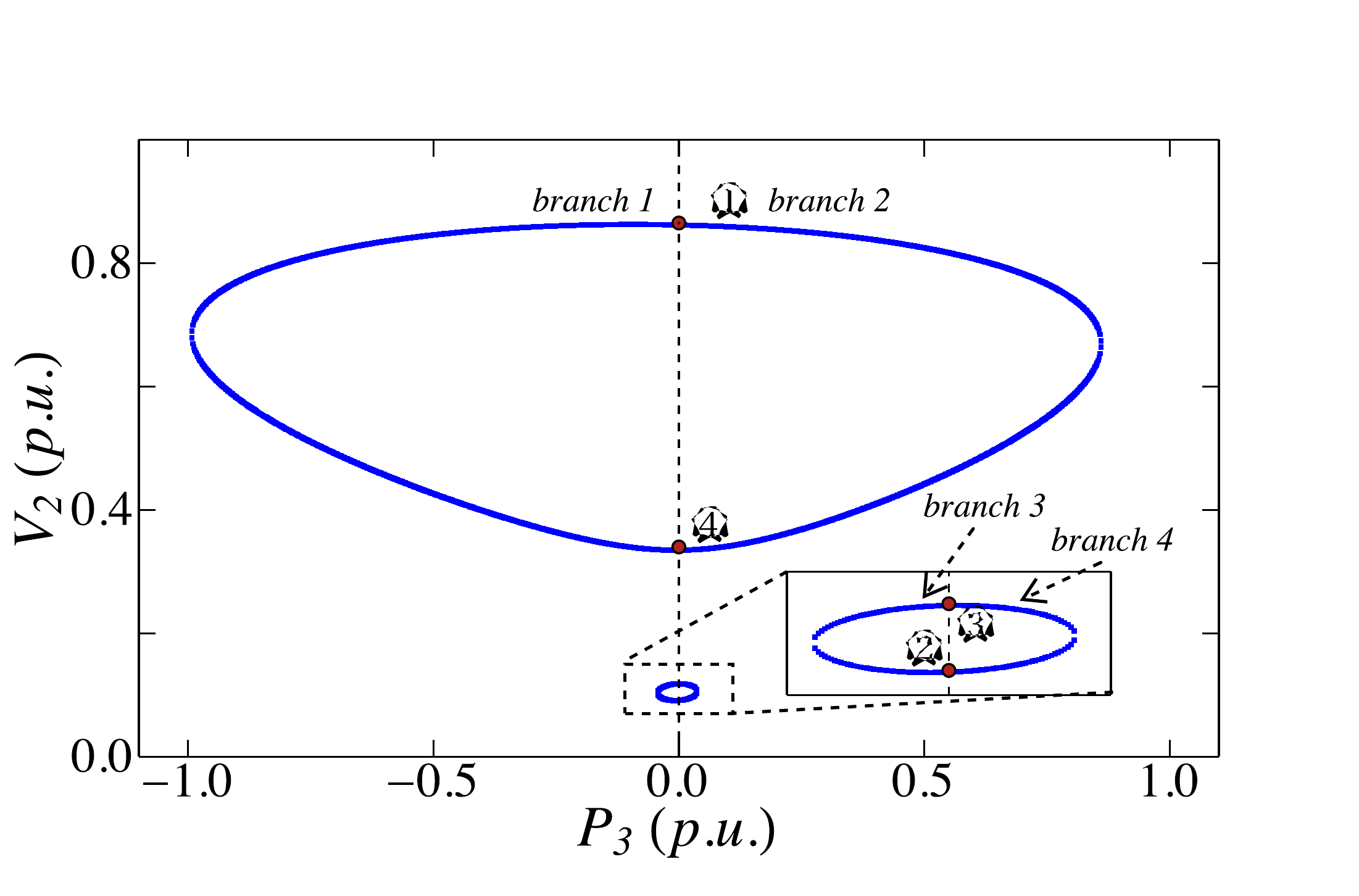}
	\caption{$P_3V_2$ curves}
    \label{fig:p3v2s}
\end{figure}

\begin{figure}[ht]
    \centering
    \includegraphics[width=1\columnwidth]{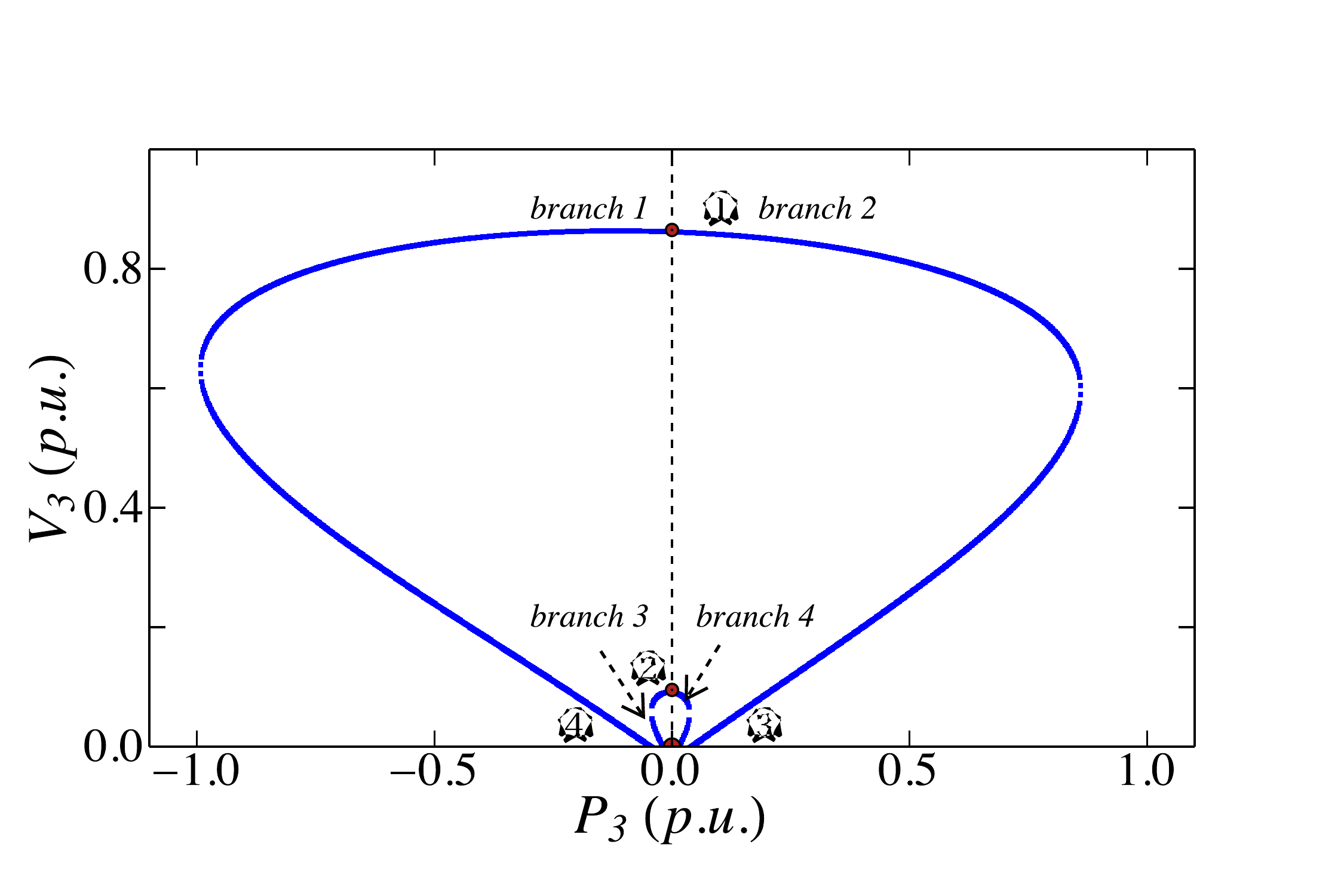}
	\caption{$P_3V_3$ curves}
    \label{fig:p3v3s}
\end{figure}

The system parameters are given as the following. Both lines $1-2$ and $3-4$ have the same impedance value  $z = 0.01 + 0.2\, p.u.$. Bus $2$ consumes $P_2 = 0.1 \,p.u.$ and $Q_2 = 0.5 \,p.u.$. Bus $3$ does not consume reactive power, i.e. $Q_3 = 0$. Active power on bus $3$ varies, with the solution branches plotted on the figures. For the $3$-bus system, there are $4$ solution branches are observed. Branches $1$ and $3$ which lie in the left half plane corresponding to bus $3$ generating active power are starting from \circled{1} and \circled{2} then ending at \circled{4} and \circled{3}, respectively. The counterparts of branches $1$ and $3$ are branches $2$ and $4$ lying in the right half plane. In Figure \ref{fig:p3v3s}, point \circled{3} and \circled{4} are congruent. Different solution branches intersect the no-loading axis $P_3=0$ at $4$ points.

At these points the power levels consumed or produced at every one of the buses are low and can be achieved either by having very low admittances (corresponding to non-ideal open circuit) or very high admittances corresponding to non-ideal short-circuit. The open and short circuits become ideal at the point where all the buses have exactly zero power consumption. The origin of this interpretation can be found in \cite{Klos1991268}.

The $4$ intersection points correspond to $4$ possible combinations of  load buses, that can be represented as either open or short circuit (s/c) statuses of the switches $s_2$ and $s_3$. The status of fictitious switches related to $4$ starting and ending points are listed in Table \ref{table:beginendpts}.

\begin{table}[ht]
    \centering
    \caption{THE CONFIGURATIONS OF THE ZERO POWER POINTS}
    \label{table:beginendpts}
    \begin{tabular}{|c |c |c |c |}
    \hline
    \textbf{Point} & \textbf{$s_2$} & \textbf{$s_3$} & \textbf{Status}\\
               \hline
    \circled{1} & $0$ & $0$ & Open both bus $2$ and $3$\\
    \circled{4} & $0$ & $1$ & Open bus $2$, s/c bus $3$\\
    \circled{2} & $1$ & $0$ & s/c bus $2$, open bus $3$\\
    \circled{3} & $1$ & $1$ & s/c both bus $2$ and $3$\\
    \hline
    \end{tabular}
\end{table}

If other active and reactive powers are non-zero such as $Q_2$, non-ideal short circuit is applied, i.e. short circuit via low impedance. Non-ideal short circuit status is corresponding to very low voltage. When the solution branch is constructed with respect to $P_3$, the status of $s_3$ will change from open to s/c; hence $g_3$ will change from $0$ to either $+\infty$ or $-\infty$ depending on the load consumes or produces power.

As observed in $P_3V_3$ curves in Figure \ref{fig:p3v3s}, all solution branches are either starting or ending at the origin of the plot which correspond to a zero power flow operating point. The key observation about the importance of the zero power flow operating point is the following. For any radial system with $n$ buses, it is possible to construct explicitly multiple different solutions of the power flow equations for the operating point where all of the power injections are zero. These solutions are constructed by either short circuiting or opening the circuit on every bus, as shown on the figure. As one can easily check the power flow in both of the configurations will be the same. So, for the system with $n-1$ $PQ$ buses we have $2^{n-1}$ possible solutions. Not all of them are different, but the important point is that many solution branches pass through this point. However, as we will argue below, most of these branches can be naturally continued far into the power generation but not in the high-power consumption regime.

Typical new solution branches in consumption regimes are tiny and close to the origin of the $PV$ plot, i.e. power transfer through the line to supply the load is limited compared to that in the production regimes. In other words, the solution branch can be usually continued into production regimes much further than in the consumption regimes. The reasons are as follows.

All the nose curves observed on the $PV$ plot are continuations of the special configurations with open or short-circuited switches. Consider a radial system where the bus $k$ is short-circuited, which means that all downstream buses with $n>k$ also have zero voltages. As continued from zero voltage level, the non-conventional solution branches will have very low voltage levels in the close neighborhood of no-loading point. Hence, they will have much lower loadability limit in the consumption region due to high currents accompanied by high power losses.

On the other hand, the solution branches can be typically indefinitely continued in the reversed power regime. The reasoning behind this conjecture is the same as for the text-book 2-bus example consisting of one slack bus and one load bus.

The load power can be expressed as $
P = V_1I - rI^2
$. If the load consumes power, i.e. $P > 0$, current $I$ will be limited when $P$ increases because $rI^2$ increases faster than $V_1I$. At the point where $I = I^{max}$ and $P = P^{max}$, the system encounters voltage collapse. However, this limit is removed in the production regimes where $P < 0$, the current $I$ can goes to infinity because power $P$ and the loss related term, $-rI^2$, have the same sign. The maximum loading in the consumption regime occurs when the impedance of load and network become comparable (equal in $2$-bus example), however this condition can never be satisfied for power flow reversal regime.

\bibliographystyle{IEEEtran}
%\bibliography{bib.bib}

\begin{thebibliography}{10}
\providecommand{\url}[1]{#1}
\csname url@samestyle\endcsname
\providecommand{\newblock}{\relax}
\providecommand{\bibinfo}[2]{#2}
\providecommand{\BIBentrySTDinterwordspacing}{\spaceskip=0pt\relax}
\providecommand{\BIBentryALTinterwordstretchfactor}{4}
\providecommand{\BIBentryALTinterwordspacing}{\spaceskip=\fontdimen2\font plus
\BIBentryALTinterwordstretchfactor\fontdimen3\font minus
  \fontdimen4\font\relax}
\providecommand{\BIBforeignlanguage}[2]{{%
\expandafter\ifx\csname l@#1\endcsname\relax
\typeout{** WARNING: IEEEtran.bst: No hyphenation pattern has been}%
\typeout{** loaded for the language `#1'. Using the pattern for}%
\typeout{** the default language instead.}%
\else
\language=\csname l@#1\endcsname
\fi
#2}}
\providecommand{\BIBdecl}{\relax}
\BIBdecl

\bibitem{Hiskens95}
I.~Hiskens, ``Analysis tools for power systems-contending with
  nonlinearities,'' \emph{Proceedings of the IEEE}, vol.~83, no.~11, pp.
  1573--1587, Nov 1995.

\bibitem{Hiskens89}
I.~Hiskens and D.~Hill, ``Energy functions, transient stability and voltage
  behaviour in power systems with nonlinear loads,'' \emph{Power Systems, IEEE
  Transactions on}, vol.~4, no.~4, pp. 1525--1533, Nov 1989.

\bibitem{dobson1992voltage}
I.~Dobson, H.~Glavitsch, C.-C. Liu, Y.~Tamura, and K.~Vu, ``Voltage collapse in
  power systems,'' \emph{Circuits and Devices Magazine, IEEE}, vol.~8, no.~3,
  pp. 40--45, 1992.

\bibitem{dobson92collapse}
I.~Dobson and L.~Lu, ``Voltage collapse precipitated by the immediate change in
  stability when generator reactive power limits are encountered,''
  \emph{Circuits and Systems I: Fundamental Theory and Applications, IEEE
  Transactions on}, vol.~39, no.~9, pp. 762--766, 1992.

\bibitem{Cutsem}
T.~Van~Cutsem and C.~Vournas, \emph{Voltage stability of electric power
  systems}.\hskip 1em plus 0.5em minus 0.4em\relax Springer, 1998, vol. 441.

\bibitem{BialekBook}
J.~Machowski, J.~Bialek, and J.~Bumby, \emph{Power system dynamics: stability
  and control}.\hskip 1em plus 0.5em minus 0.4em\relax John Wiley \& Sons,
  2011.

\bibitem{Kundur}
P.~Kundur, \emph{{Power System Stability and Control}}, New York, 1994.

\bibitem{Thorp}
W.~Ma and J.~Thorp, ``An efficient algorithm to locate all the load flow
  solutions,'' \emph{Power Systems, IEEE Transactions on}, vol.~8, no.~3, pp.
  1077--1083, Aug 1993.

\bibitem{Chiang90}
H.-D. Chiang and M.~Baran, ``On the existence and uniqueness of load flow
  solution for radial distribution power networks,'' \emph{Circuits and
  Systems, IEEE Transactions on}, vol.~37, no.~3, pp. 410--416, 1990.

\bibitem{thorp1986reactive}
J.~Thorp, D.~Schulz, and M.~Ili{\'c}-Spong, ``Reactive power-voltage problem:
  conditions for the existence of solution and localized disturbance
  propagation,'' \emph{International Journal of Electrical Power \& Energy
  Systems}, vol.~8, no.~2, pp. 66--74, 1986.

\bibitem{ilic1992network}
M.~Ilic, ``Network theoretic conditions for existence and uniqueness of steady
  state solutions to electric power circuits,'' in \emph{Circuits and Systems,
  1992. ISCAS'92. Proceedings., 1992 IEEE International Symposium on},
  vol.~6.\hskip 1em plus 0.5em minus 0.4em\relax IEEE, 1992, pp. 2821--2828.

\bibitem{wang2012distflow}
D.~Wang, K.~Turitsyn, and M.~Chertkov, ``Distflow ode: Modeling, analyzing and
  controlling long distribution feeder,'' in \emph{Decision and Control (CDC),
  2012 IEEE 51st Annual Conference on}.\hskip 1em plus 0.5em minus 0.4em\relax
  IEEE, 2012, pp. 5613--5618.

\bibitem{Overbye1994}
T.~Overbye, ``Effects of load modelling on analysis of power system voltage
  stability,'' \emph{International Journal of Electrical Power \& Energy
  Systems}, vol.~16, no.~5, pp. 329--338, 1994.

\bibitem{Nose_Taranto}
S.~Corsi and G.~Taranto, ``Voltage instability - the different shapes of the
  "nose",'' in \emph{Bulk Power System Dynamics and Control - VII. Revitalizing
  Operational Reliability, 2007 iREP Symposium}, Aug 2007, pp. 1--16.

\bibitem{Venkatasubramanian1992}
V.~Venkatasubramanian, H.~Schattler, and J.~Zaborszky, ``Voltage dynamics:
  study of a generator with voltage control, transmission, and matched mw
  load,'' \emph{Automatic Control, IEEE Transactions on}, vol.~37, no.~11, pp.
  1717--1733, Nov 1992.

\bibitem{hill1994stability}
D.~Hill, I.~Hiskens, and D.~Popovic, ``Stability analysis of power system loads
  with recovery dynamics,'' \emph{International Journal of Electrical Power \&
  Energy Systems}, vol.~16, no.~4, pp. 277--286, 1994.

\bibitem{Hiskens2010}
R.~Bravo, R.~Yinger, D.~Chassin, H.~Huang, N.~Lu, I.~Hiskens, and
  G.~Venkataramanan, ``Final project report load modeling transmission
  research,'' \emph{Lawrence Berkeley National Laboratory (LBNL)}, 2010.

\bibitem{Hiskens2005}
V.~Donde and I.~Hiskens, ``Dynamic performance assessment: grazing and related
  phenomena,'' \emph{Power Systems, IEEE Transactions on}, vol.~20, no.~4, pp.
  1967--1975, Nov 2005.

\bibitem{Dobson2012}
H.~Wu and I.~Dobson, ``Cascading stall of many induction motors in a simple
  system,'' \emph{Power Systems, IEEE Transactions on}, vol.~27, no.~4, pp.
  2116--2126, Nov 2012.

\bibitem{turitsyn2011options}
K.~Turitsyn, P.~Sulc, S.~Backhaus, and M.~Chertkov, ``Options for control of
  reactive power by distributed photovoltaic generators,'' \emph{Proceedings of
  the IEEE}, vol.~99, no.~6, pp. 1063--1073, 2011.

\bibitem{Karrison94}
D.~Karlsson and D.~Hill, ``Modelling and identification of nonlinear dynamic
  loads in power systems,'' \emph{Power Systems, IEEE Transactions on}, vol.~9,
  no.~1, pp. 157--166, Feb 1994.

\bibitem{Hill93}
D.~J. Hill, M.~Pal, X.~Wilsun, Y.~Mansour, C.~Nwankpa, L.~Xu, and R.~Fischl,
  ``Nonlinear dynamic load models with recovery for voltage stability studies.
  discussion. authors' response,'' \emph{IEEE Transactions on Power Systems},
  vol.~8, no.~1, pp. 166--176, 1993.

\bibitem{Erlich2006}
I.~Erlich, K.~Rensch, and F.~Shewarega, ``Impact of large wind power generation
  on frequency stability,'' in \emph{Power Engineering Society General Meeting,
  2006. IEEE}.\hskip 1em plus 0.5em minus 0.4em\relax IEEE, 2006, pp. 8--pp.

\bibitem{Tamimi2013}
B.~Tamimi, C.~Canizares, and K.~Bhattacharya, ``System stability impact of
  large-scale and distributed solar photovoltaic generation: The case of
  ontario, canada,'' \emph{Sustainable Energy, IEEE Transactions on}, vol.~4,
  no.~3, pp. 680--688, July 2013.

\bibitem{IEEEtestcase}
\BIBentryALTinterwordspacing
``Distribution test feeders.'' [Online]. Available:
  \url{http://ewh.ieee.org/soc/pes/dsacom/testfeeders/index.html}
\BIBentrySTDinterwordspacing

\bibitem{FIDVRwhitepaper}
W.~E. C.~C. Modeling and V.~W. Group, ``Modeling and studying fidvr events,''
  October 20, 2011.

\bibitem{Wangthesis}
J.~Wang, ``A planning scheme for penetrating embedded generation in power
  distribution grids,'' Ph.D. dissertation, Massachusetts Institute of
  Technology, Massachusetts, Cambridge, MA, 2013.

\bibitem{milano2011hybrid}
F.~Milano, ``Hybrid control model of under load tap changers,'' \emph{Power
  Delivery, IEEE Transactions on}, vol.~26, no.~4, pp. 2837--2844, 2011.

\bibitem{Klos1991268}
{A. Klos and J. Wojcicka}, ``Physical aspects of the nonuniqueness of load flow
  solutions,'' \emph{International Journal of Electrical Power \& Energy
  Systems}, vol.~13, no.~5, pp. 268 -- 276, 1991.

\end{thebibliography}

\begin{IEEEbiography}[{\includegraphics[width=1in,height=1.25in,clip,keepaspectratio]{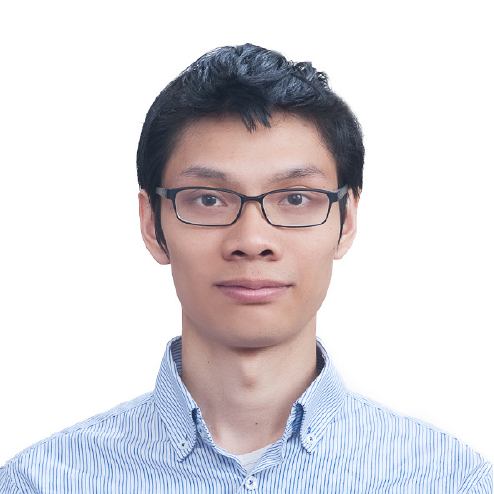}}]{Hung D. Nguyen} (S`12) was born in Vietnam, in 1986. He received the B.E. degree in electrical engineering from Hanoi University of Technology, Vietnam, in 2009, and the M.S. degree in electrical engineering from Seoul National University, Korea, in 2013. He is pursuing a Ph.D. degree in the Department of Mechanical Engineering at Massachusetts Institute of Technology (MIT). His current research interests include power system operation and control; the nonlinearity, dynamics and stability of large scale power systems; DSA/EMS and smart grids.
\end{IEEEbiography}
%\vspace{-10 mm}
\newpage
\begin{IEEEbiography}[{\includegraphics[width=1in,height=1.25in,clip,keepaspectratio]{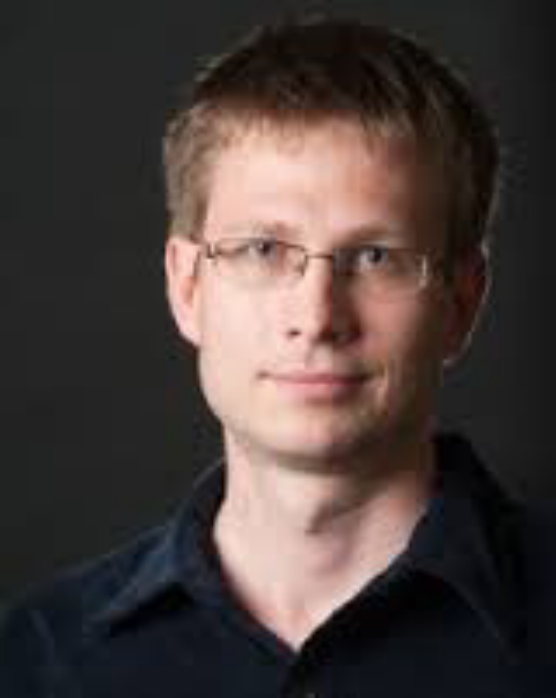}}]{Konstantin Turitsyn} (M`09) received the M.Sc. degree in physics from Moscow Institute of Physics and Technology and the Ph.D. degree in physics from Landau Institute for Theoretical Physics, Moscow, in 2007.  Currently, he is an Assistant Professor at the Mechanical Engineering Department of Massachusetts Institute of Technology (MIT), Cambridge. Before joining MIT, he held the position of Oppenheimer fellow at Los Alamos National Laboratory, and Kadanoff–Rice Postdoctoral Scholar at University of Chicago. His research interests encompass a broad range of problems involving nonlinear and stochastic dynamics of complex systems. Specific interests in energy related fields include stability and security assessment, integration of distributed and renewable generation.
\end{IEEEbiography}

\end{document}